\documentstyle[12pt]{article}
\def\mineappendix{
        \setcounter{section}{1}
        \setcounter{subsection}{0}
        \def\thesection{\Alph{section}}
        \def\sectionap{\@startsection  {section}{1}{\z@}
                        {-3.5ex plus-1ex minus-.2ex} {0ex plus.2ex}
                        {\reset@font\Large\bf  Appendix:  \, }
                        }
        }
\makeatother
\def\Proclaim #1. #2\par{\bigbreak\noindent{\sc#1.\enspace}{\it#2}\par}

\font\Bbbfont=msbm10
\newfam\msbfam
\textfont\msbfam=\Bbbfont
\def\Bbb#1{{\fam\msbfam\relax#1}}

\newcommand{\eqref}[1]{equation~(\ref{#1})}
\newcommand{\Eqref}[1]{Equation~(\ref{#1})}

\newcommand{\gwii}[1]{\left< \hspace{-2pt} \left< \, #1 \,
        \right>  \hspace{-2pt} \right>_{0}}

\newcommand{\gwiione}[1]{\left< \hspace{-2pt} \left< \, #1 \,
        \right> \hspace{-2pt} \right>_{1}}

\newcommand{\gwiitwo}[1]{\left< \hspace{-2pt} \left< \, #1 \,
        \right> \hspace{-2pt} \right>_{2}}

\newcommand{\gwig}[1]{\left< \, #1 \, \right>_{g}}
\newcommand{\gwiig}[1]{\left< \hspace{-2pt} \left< \, #1 \,
    \right> \hspace{-2pt} \right>_{g}}

\newcommand{\grav}[2]{\tau_{#1}(\gamma_{#2})}
\newcommand{\grava}[1]{\tau_{#1}(\gamma_{\alpha})}

\newcommand{\gravb}[1]{\tau_{#1}(\gamma_{\beta})}

\newcommand{\ga}{\gamma_{\alpha}}
\newcommand{\gua}{\gamma^{\alpha}}
\newcommand{\gb}{\gamma_{\beta}}
\newcommand{\gub}{\gamma^{\beta}}

\newcommand{\vs}{{\cal S}}
\newcommand{\vx}{{\cal X}}
\newcommand{\vd}{{\cal D}}
\newcommand{\ve}{{\cal E}}
\newcommand{\vf}{{\cal F}}
\newcommand{\vw}{{\cal W}}
\newcommand{\vv}{{\cal V}}

\newcommand{\vl}{{\cal L}}
\newcommand{\vg}{{\cal G}}

\newcommand{\bvs}{{\, \overline{\cal S}\,}}
\newcommand{\bvx}{{\, \overline{\cal X}\,}}

\newcommand{\bvw}{{\, \overline{\cal W} \,}}

\newcommand{\bvl}{{\, \overline{\cal L} \,}}

\newcommand{\qp} {\circ}

\newtheorem{lem}{Lemma}[section]
\newtheorem{cor}[lem]{Corollary}
\newtheorem{thm}[lem]{Theorem}

\title{Genus-2 Gromov-Witten invariants for manifolds with semisimple
         quantum cohomology}
\author{Xiaobo Liu \thanks{Research partially supported by
            Alfred P. Sloan Research Fellowship and National Science
        Foundation research grant}}
\date{}

\begin{document}
\maketitle

In \cite{L2}, the author studied universal equations for genus-2
Gromov-Witten invariants given in \cite{Ge1} and \cite{BP} using
quantum product on the big phase space. Among other results, the
author proved that for manifolds with semisimple quantum
cohomology, the generating function for
genus-2 Gromov-Witten invariants, denoted by  $F_{2}$, 
is uniquely determined by known genus-2 universal equations.
Moreover, an explicit formula for $F_{2}$ was given in terms of
genus-0 and genus-1 invariants. However, the formula given in
\cite{L2} is very complicated to work with. In this paper, we
will give a much simpler formula using idempotents of the quantum
product on the big phase space, and then use it to prove the
genus-2 Virasoro conjecture  for manifolds with semisimple quantum
cohomology (cf. \cite{EHX} and \cite{CK}).

Properties of idempotents on the big phase space were studied in
\cite{L4}. Let $M$ be a compact symplectic manifold. In
Gromov-Witten theory, the space $H^{*}(M; {\Bbb C})$ is called the
{\it small phase space}. A product of infinitely many copies of
the small phase space is called the {\it big phase space}. The
generating functions for Gromov-Witten invariants are formal
power series on the big phase space. Let $N$ be the dimension of
$H^{*}(M; {\Bbb C})$. If the quantum cohomology of $M$ is
semisimple, there exist vector fields $\ve_{i}$, $i=1, \ldots, N$,
 on the big phase space
such that $\ve_{i} \qp \ve_{j} = \delta_{i j} \ve_{i}$
for all $i$ and $j$, where ``$\qp$'' stands for the quantum product
 (see \eqref{eqn:QP}).  These vector fields
are called {\it idempotents}. Let $u_{i}$, $i=1, \ldots, N$, be
the eigenvalues of the quantum multiplication by the Euler vector
field (see \eqref{eqn:X}). The first main result of this paper is
the following
\begin{thm} \label{thm:F2}
If the quantum cohomology of the underlying manifold is semisimple,
the genus-2 generating function $F_{2}$ is given by
\[ F_{2} = \frac{1}{2} A_{1}(\tau_{-}(\vs)) +
    \frac{1}{3} A_{1}(\tau_{-}^{2}(\vl_{0}))
    - \frac{1}{6} \sum_{i=1}^{N} u_{i} B(\ve_{i}, \ve_{i}, \ve_{i}).\]
\end{thm}
In this formula, $A_{1}$ and $B$ are tensors which only depend on
genus-0 and genus-1 data. Precise formulas for $A_{1}$ and
$B(\ve_{i}, \ve_{i}, \ve_{i})$ are given in \eqref{eqn:A1}
and \eqref{eqn:Bi} respectively.
The tensor $A_{1}$ comes from the genus-0
and genus-1 part of the genus-2 topological recursion relation
derived from Mumford's relation (cf. \cite{Ge1}) and tensor $B$ comes
from the
genus-0 and genus-1 part of an equation due to Belorousski and
Pandharipande (cf. \cite{BP}). $\vs$ is the string vector field
(see \eqref{eqn:S}), and $\vl_{n}$ is the $n$-th Virasoro vector
field (see \eqref{eqn:VirVect}). The operator $\tau_{-}$ lower the
level of descendants of vector fields by 1.
We will also give a formula for $F_{2}$
which only involves genus-0 data in Theorem~\ref{thm:F2rot}.

The Virasoro conjecture predicts that the generating functions of
Gromov-Witten invariants of smooth projective varieties are
annihilated by an infinite sequence of differential operators
which form a half branch of the Virasoro algebra. This conjecture
was proposed by Eguchi-Hori-Xiong and Katz (cf. \cite{EHX},
\cite{CK}). It is a natural generalization of Witten's KdV
conjecture (cf. \cite{W1} \cite{W2}) which was proved by
Kontsevich (cf. \cite{Kon}). For manifolds with semisimple quantum
cohomology, the Virasoro conjecture completely determines the
higher genus Gromov-Witten invariants in terms of genus-0
invariants (cf. \cite{DZ3}). In \cite{LT}, Tian and the author
proved the genus-0 Virasoro conjecture for all compact symplectic
manifolds (see also \cite{DZ2}, \cite{Ge2}, \cite{L3},
\cite{Gi3}). The genus-1 Virasoro conjecture for manifolds with
semisimple quantum cohomology was proved by Dubrovin and Zhang
\cite{DZ2} (see also \cite{L1} and \cite{L4}). In \cite{L1} and
\cite{L2}, the author also proved that the genus-1 and genus-2
Virasoro conjecture for all smooth projective varieties can be
reduced to an $SL(2)$ symmetry of Gromov-Witten invariants. The
second main result of this paper is the following
\begin{thm} \label{thm:g2VirSS}
For smooth projective varieties with semisimple quantum cohomology,
the genus-2 Virasoro conjecture is true.
\end{thm}
In \cite{Gi1}, Givental conjectured a formula for higher genus
Gromov-Witten potential for manifolds with semisimple quantum
cohomology. His formula satisfies the Virasoro constraints (cf.
\cite{Gi2}). Since in
the semisimple case, Virasoro constraints uniquely determine the
higher genus Gromov-Witten potential (cf. \cite{DZ3}),
Theorem~\ref{thm:g2VirSS} implies
that Givental's conjectural formula is correct in the genus-2
case. We also note that the method used in this paper should apply to 
higher genus case once the corresponding universal equations are obtained.

Part of the work in this paper was done when the author visited
IPAM at Los Angeles, MSRI at Berkeley, and IHES in France. The
author would like to thank these institutes for hospitality.

\section{Quantum product and idempotents}

We first review properties of quantum product and idempotents on
the big phase space which will be used in this paper. The proofs
for these properties can be found in  \cite{L2} and \cite{L4}.

Let $M$ be a compact symplectic manifold. For simplicity, we assume
$H^{\rm odd}(M; {\Bbb C}) = 0$.
The {\it big phase space} is by definition the product of infinite copies of
$H^{*}(M; {\Bbb C})$, i.e.
\[  P := \prod_{n=0}^{\infty} H^{*}(M; {\Bbb C}). \]
Fix a basis $\{ \gamma_{1}, \ldots, \gamma_{N} \}$ of
$H^{*}(M; {\Bbb C})$ with $\gamma_{1} = 1$ being the identity of the ordinary
cohomology ring of $M$. Then we denote the corresponding basis for
the $n$-th copy of $H^{*}(M; {\Bbb C})$ in $P$ by
$\{\tau_{n}(\gamma_{1}), \ldots, \tau_{n}(\gamma_{N}) \}$.
We call $\grava{n}$ a {\it descendant} of $\gamma_{\alpha}$ with descendant
level $n$.
We can think of $P$ as an infinite dimensional vector space with basis
$\{ \grava{n} \mid 1 \leq \alpha \leq N, \, \, \, n \in {\Bbb Z}_{\geq 0} \}$
where ${\Bbb Z}_{\geq 0} = \{ n \in {\Bbb Z} \mid n \geq 0\}$.
Let
$(t_{n}^{\alpha} \mid 1 \leq \alpha \leq N, \, \, \, n \in {\Bbb Z}_{\geq 0})$
be the corresponding coordinate system on $P$.
For convenience, we identify $\grava{n}$ with the coordinate vector field
$\frac{\partial}{\partial t_{n}^{\alpha}}$ on $P$ for $n \geq 0$.
If $n<0$, $\grava{n}$ is understood as the $0$ vector field.
We also abbreviate $\grava{0}$ as $\gamma_{\alpha}$.
Any vector field of the form $\sum_{\alpha} f_{\alpha} \ga$,
where $f_{\alpha}$
are functions on the big phase space, is called a {\it primary vector field}.
We use $\tau_{+}$ and $\tau_{-}$ to denote the operator which shift the level
of descendants, i.e.
\[ \tau_{\pm} \left(\sum_{n, \alpha} f_{n, \alpha} \grava{n}\right)
    = \sum_{n, \alpha} f_{n, \alpha} \grava{n \pm 1} \]
where $f_{n, \alpha}$ are functions on the big phase space.

We will use the following {\it conventions for notations}:  All
summations are over the entire meaningful ranges of the indices
unless otherwise indicated. Let
\[ \eta_{\alpha \beta} = \int_{M} \gamma_{\alpha} \cup
    \gamma_{\beta}
\]
 be the intersection form on $H^{*}(M, {\Bbb C})$.
We will use $\eta = (\eta_{\alpha \beta})$ and $\eta^{-1} =
(\eta^{\alpha \beta})$ to lower and raise indices.
For example,
\[ \gua :=  \eta^{\alpha \beta} \gb.\]
Here we are using the summation convention that repeated
indices (in this formula, $\beta$) should be summed
over their entire ranges.

Let
\[ \gwig{\grav{n_{1}}{\alpha_{1}} \, \grav{n_{2}}{\alpha_{2}} \,
    \ldots \, \grav{n_{k}}{\alpha_{k}}} \]
be the genus-$g$ descendant Gromov-Witten invariant associated
to $\gamma_{\alpha_{1}}, \ldots, \gamma_{\alpha_{k}}$ and nonnegative
integers $n_{1}, \ldots, n_{k}$
(cf. \cite{W1}, \cite{RT}, \cite{LiT}).
The genus-$g$
generating function is defined to be
\[ F_{g} =  \sum_{k \geq 0} \frac{1}{k!}
         \sum_{ \begin{array}{c}
        {\scriptstyle \alpha_{1}, \ldots, \alpha_{k}} \\
                {\scriptstyle  n_{1}, \ldots, n_{k}}
                \end{array}}
                t^{\alpha_{1}}_{n_{1}} \cdots t^{\alpha_{k}}_{n_{k}}
    \gwig{\grav{n_{1}}{\alpha_{1}} \, \grav{n_{2}}{\alpha_{2}} \,
        \ldots \, \grav{n_{k}}{\alpha_{k}}}. \]
This function is understood as a formal power series of
$t_{n}^{\alpha}$.

Introduce
a $k$-tensor
 $\left< \left< \right. \right. \underbrace{\cdot \cdots \cdot}_{k}
        \left. \left. \right> \right> $
defined by
\[ \gwiig{{\cal W}_{1} {\cal W}_{2} \cdots {\cal W}_{k}} \, \,
         := \sum_{m_{1}, \alpha_{1}, \ldots, m_{k}, \alpha_{k}}
                f^{1}_{m_{1}, \alpha_{1}} \cdots f^{k}_{m_{k}, \alpha_{k}}
        \, \, \, \frac{\partial^{k}}{\partial t^{\alpha_{1}}_{m_{1}}
            \partial t^{\alpha_{2}}_{m_{k}} \cdots
            \partial t^{\alpha_{k}}_{m_{k}}} F_{g},
 \]
for  vector fields ${\cal W}_{i} = \sum_{m, \alpha} f^{i}_{m,
\alpha} \, \frac{\partial}{\partial t_{m}^{\alpha}}$ where
$f^{i}_{m, \alpha}$ are functions on the big phase space. We can
also view this tensor as the $k$-th covariant derivative of
$F_{g}$ with respect to the trivial connection on $P$. This tensor
is called the {\it $k$-point (correlation) function}. For any
vector fields $\vw_{1}$ and $\vw_{2}$ on the big phase space, the
{\it quantum product}
 of $\vw_{1}$ and $\vw_{2}$ is defined by
\begin{equation} \label{eqn:QP}
 \vw_{1} \qp \vw_{2} := \gwii{\vw_{1} \, \vw_{2} \, \gua} \ga.
\end{equation}
This is a commutative and associative product. But it does not have
an identity. For any vector field $\vw$ and integer $k \geq 1$,
$\vw^{k}$ is understood as the $k$-th power of $\vw$ with respect to this
product.

Let
\begin{equation} \label{eqn:X}
{\vx} := - \sum_{m, \alpha} \left(m + b_{\alpha} - b_{1} -1
                        \right)\tilde{t}^{\alpha}_{m} \, \grava{m}
        - \sum_{m, \alpha, \beta}
        {\cal C}_{\alpha}^{\beta}\tilde{t}^{\alpha}_{m} \,
        \gravb{m-1}
\end{equation}
be the Euler vector field on the big phase space $P$,
where 
$\tilde{t}^{\alpha}_{m} = t^{\alpha}_{m} - \delta_{m,1} \delta_{\alpha, 1}$,
\[ b_{\alpha} = \frac{1}{2}({\rm dimension \,\,\, of \,\,\,} \ga)
   - \frac{1}{4} ({\rm real \,\,\, dimension \,\,\, of \,\,\,} M)
   + \frac{1}{2}\]
 and the matrix ${\cal C} = ( {\cal C}_{\alpha}^{\beta})$
is defined by $ c_{1}(V) \cup \ga = {\cal C}_{\alpha}^{\beta} \,
\gb$. For smooth projective varieties, the dimension of $\ga$
should be replaced by twice of the holomorphic dimension of $\ga$
in the definition of $b_{\alpha}$.

The quantum multiplication by $\vx$ is an endomorphism on the space
of primary vector fields on $P$.
If this endomorphism has distinct eigenvalues at generic points,
we call $P$ {\it semisimple}. In this case, let $\ve_{1}, \ldots, \ve_{N}$
be the eigenvectors with corresponding eigenvalues $u_{1}, \ldots, u_{N}$,
i.e.
\[ \vx \qp \ve_{i} = u_{i} \ve_{i} \]
for each $i = 1, \cdots, N$. $\ve_{i}$ is considered as a vector field on
$P$, and $u_{i}$ is considered as a  function on $P$.
They satisfy the following properties:
\[ \ve_{i} \qp \ve_{j} = \delta_{ij} \, \ve_{i}, \hspace{20pt}
      [ \ve_{i}, \,  \ve_{j}] = 0, \hspace{20pt}
     \ve_{i} \, u_{j} = \delta_{ij}
\]
for any $i$ and $j$.
We call $\{ \ve_{1}, \ldots, \ve_{N} \}$  {\it idempotents} on the
big phase space. When restricted to the small phase space,
they coincide with the coordinate vector fields of the canonical
coordinate system of semisimple Frobenius manifolds (cf. \cite{D}).

Let
\begin{equation} \label{eqn:S}
 {\vs} := - \sum_{m, \alpha} \tilde{t}^{\alpha}_{m}
        \grav{m-1}{\alpha}
\end{equation}
be the {\it string vector field} on $P$. We define
\[ \bvw = \vw \qp \vs \]
for any vector field $\vw$ on $P$.
The vector field $\bvs$ is the identity for the quantum product 
when restricted 
to the space of primary vector fields.
We have
\begin{equation} \label{eqn:bvs}
     \bvs = \sum_{i=1}^{N} \ve_{i}.
\end{equation}
and
\begin{equation} \label{eqn:bvx}
 \bvx^{k} = \sum_{i=1}^{N} u_{i}^{k} \, \ve_{i}.
\end{equation}
for $k \geq 1$.

For any vector fields $\vw$ and $\vv$ on the big phase space, define
\begin{equation} \label{eqn:inner}
<\vw, \vv> := \gwii{\vs \, \vw \, \vv}.
\end{equation}
This bilinear form generalizes the Poincare
metric on the small phase space. It is nondegenerate only
when restricted to the space of primary vector fields.
It is also compatible with quantum product in the
following sense:
\begin{equation}
 <(\vw_{1} \qp \vw_{2}), \, \vw_{3}>
          = < \vw_{2}, \, (\vw_{1} \qp \vw_{3})>, \label{eqn:Frob}
\end{equation}
 and
\[ <\vw_{1}, \vw_{2}> = <\bvw_{1}, \bvw_{2}> \]
for any vector fields $\vw_{i}$.
We have $<\ve_{i}, \ve_{j}> = 0$ if $i \neq j$ and the functions
\[ g_{i} := <\ve_{i}, \, \ve_{i}> \]
are non-zero in the region where idempotents are well defined.
Any primary vector field $\vw$ has the decomposition
\[ \vw = \sum_{i=1}^{N} \frac{<\vw, \, \ve_{i}>}{g_{i}} \, \ve_{i}. \]

Let $\nabla$ be the covariant derivative on $P$ of the trivial flat connection
with respect to the standard coordinates $\{t_{n}^{\alpha}\}$.
The compatabilities of this connection with the quantum product and
the bilinear form are given by
\begin{equation}
  \nabla_{\vw_{1}} (\vw_{2} \qp \vw_{3})
= (\nabla_{\vw_{1}} \vw_{2}) \qp \vw_{3}
    + \vw_{2} \qp (\nabla_{\vw_{1}} \vw_{3})
    + \gwii{\vw_{1} \, \vw_{2} \, \vw_{3} \, \gua} \ga 
\label{eqn:connqp}
\end{equation}
and
\begin{eqnarray}
 \vw_{1} <\vw_{2}, \vw_{3}>
&=& < \left\{ \nabla_{\vw_{1}} \vw_{2}
        + \vw_{1} \qp \tau_{-}(\vw_{2}) \right\}, \vw_{3}>  \nonumber \\
&&    + <\vw_{2},  \left\{ \nabla_{\vw_{1}} \vw_{3}
                    + \vw_{1} \qp \tau_{-}(\vw_{3}) \right\}>
                    \label{eqn:innermetric}
\end{eqnarray}
for any vector fields $\vw_{i}$. \Eqref{eqn:innermetric} suggests that
the modified connection $\tilde{\nabla}$ defined by
\[ \tilde{\nabla}_{\vw_{1}} \vw_{2} :=
	\nabla_{\vw_{1}} \vw_{2}
        + \vw_{1} \qp \tau_{-}(\vw_{2}) \]
is compatible with the bilinear form $< \cdot, \, \cdot>$. Moreover
\eqref{eqn:connqp} implies that the family of connections 
$\tilde{\nabla}^{z}$ defined by
\[ \tilde{\nabla}^{z}_{\vw_{1}} \vw_{2} :=
	\nabla_{\vw_{1}} \vw_{2}
        + z \, \vw_{1} \qp \tau_{-}(\vw_{2}) \]
are flat for all $z$, where $z$ is an arbitrary parameter.

Covariant derivatives of idempotents
are given by
\[ \nabla_{_{\vw}} \ve_{i}
= -2 \gwii{\vw \, \ve_{i} \, \ve_{i} \, \gua} \ga \qp \ve_{i}
        + \gwii{\vw \, \ve_{i} \, \ve_{i} \, \gua} \ga \]
for any vector field $\vw$. In particular,
\[ \nabla_{\ve_{j}} \ve_{i} = \delta_{ij} \vf_{j} -
     \vf_{j} \qp \ve_{i} -   \vf_{i} \qp \ve_{j} \]
where
\[ \vf_{j} := \gwii{\ve_{j} \, \ve_{j} \, \ve_{j} \, \gua} \ga \]
for each $j = 1, \cdots, N$.
Vector fields $\vf_{j}$ are also related to the string vector
field by the formula
\begin{equation} \label{eqn:tau-S}
\overline{\tau_{-}(\vs)} = - \sum_{i=1}^{N} \vf_{i}.
\end{equation}

For any vector field $\vw$, define
\[ T(\vw) := \tau_{+}(\vw) - \vs \qp \tau_{+}(\vw). \]
The operator $T$ was introduced in \cite{L2} to simplify
topological recursion relations for Gromov-Witten invariants. It
corresponds to the $\psi$ classes in the relations in the
tautological ring of moduli space of stable curves. In some
sense, repeatedly applying $T$ to a vector field will trivialize
its action on genus-$g$ generating functions. Here are some basic
properties of $T$: For any vector fields $\vw_{i}$,
\begin{eqnarray*}
&(i)& T(\vw_{1}) \qp \vw_{2} = 0, \\
&(ii)&   \gwii{T(\vw_{1}) \, \vw_{2} \, \vw_{3} \, \vw_{4}}
    = \gwii{(\vw_{1} \qp \vw_{2}) \, \vw_{3} \, \vw_{4}} \\
&(iii)&
    \nabla_{\vw_{1}} \,\, T(\vw_{2}) = T \left(\nabla_{\vw_{1}} \vw_{2}
        \right) - \vw_{1} \qp \vw_{2} \\
&(iv)& T(\vw) \, u_{i} = 0, \\
&(v)& <T(\vw_{1}), \vw_{2}>  = 0, \\
&(vi)& \nabla_{_{T(\vw)}} \ve_{i} = - \vw \qp \ve_{i}, \\
&(vii)& [T(\vw), \, \ve_{i}] = - T(\nabla_{\ve_{i}} \vw).
\end{eqnarray*}
Note that $\{T^{k}(\ve_{i}) \mid i=1, \ldots, N, k \geq 0 \}$
gives a frame for the tangent bundle of the big phase space. This
frame is not commutative due to property (vii). Any vector field
$\vw$ has the following decomposition
\begin{equation} \label{eqn:WTW}
\vw = T^{k} (\tau_{-}^{k}(\vw)) + \sum_{i=0}^{k-1} T^{i}
(\overline{\tau_{-}^{i}(\vw)})
\end{equation}
where $k$ is any positive integer (cf \cite[Equation (26)]{L2}).
This decomposition
is very useful when applying topological recursion relations.
In particular, we will frequently use the decomposition
\[ \vw = \bvw + T(\tau_{-}(\vw)) \]
and call this the {\it standard decomposition} of $\vw$.
For example,  using this decomposition, we see
\[ \vw \, u_{i} = \bvw \, u_{i} \]
for any vector field $\vw$ on the big phase space.

For any vector field
$\vw = \sum_{n, \alpha} f_{n, \alpha} \tau_{n}(\gamma_{\alpha})$,
define
\[ \vg * \vw := \sum_{n, \alpha}
    (n + b_{\alpha}) f_{n, \alpha} \tau_{n}(\gamma_{\alpha}). \]
This operator was used in \cite{L2} to give a recursive description
for the Virasoro vector fields.
On the space of primary vector fields, the operator $\vg *$
has the following property:
\begin{equation} \label{eqn:G*inner}
 <\vg * \vw, \, \vv> + <\vw, \, \vg * \vv> = <\vw, \, \vv>
\end{equation}
for all primary vector fields $\vv$ and $\vw$.
Moreover, for any $i$,
\begin{equation} \label{eqn:g*ei}
 \vg * \ve_{i} = \frac{1}{2} \, \ve_{i}
        - u_{i} \vf_{i} + \vx \qp \vf_{i}.
\end{equation}

\section{Universal equations in genus 2}

We will need two genus-2 universal equations. The first one is the
genus-2 topological recursion relation derived from Mumford's relation
 (cf. \cite{Ge1}):
For any vector field $\vw$,
\begin{equation} \label{eqn:TRR1}
\gwiitwo{T^{2}({\cal W})} = A_{1}(\vw) \hspace{200pt}
\end{equation}
where
\begin{eqnarray}
A_{1}(\vw) &:=&
    \frac{7}{10} \gwiione{\ga} \gwiione{\{\gua \qp \vw\}}
 + \frac{1}{10} \gwiione{\ga \, \{\gua \qp \vw\}}
     - \frac{1}{240} \gwiione{\vw \, \{\ga \qp \gua\}} 
	\nonumber \\
&&    + \frac{13}{240} \gwii{\vw \, \ga \, \gua \, \gub}
        \gwiione{\gb}
  +   \frac{1}{960} \gwii{ \vw \, \gua \, \ga \, \gub \, \gb}.
	\label{eqn:A1}
\end{eqnarray}
Another genus-2 equation is the following (cf. \cite{BP}): For any
vector fields $\vw_{i}$,
\begin{eqnarray}
&&  2 \gwiitwo{\{{\cal W}_{1} \qp {\cal W}_{2} \qp {\cal W}_{3} \}}
    - 2 \gwii{{\cal W}_{1} \, {\cal W}_{2} \, {\cal W}_{3} \, \gua}
    \gwiitwo{T(\ga)} \nonumber \\
&& + \frac{1}{2} \sum_{\sigma \in S_{3}} \left\{
    \gwiitwo{{\cal W}_{\sigma(1)} \, T({\cal W}_{\sigma(2)}
        \qp {\cal W}_{\sigma(3)})}
    - \gwiitwo{ T({\cal W}_{\sigma(1)}) \, \{{\cal W}_{\sigma(2)}
        \qp {\cal W}_{\sigma(3)}\}} \right\} \nonumber \\
& = &  B(\vw_{1}, \vw_{2}, \vw_{3}), \label{eqn:BP}
\end{eqnarray}
where $B$ is a symmetric 3-tensor which only
depends on genus 0 and genus 1 data.
The precise definition for $B$ is very complicated (see \cite[Section 2]{L2}).
In this paper, we only need the special case $B(\ve_{i}, \ve_{i}, \ve_{i})$.
We will give the definition for this function in \eqref{eqn:Bi}.

We will also need the Virasoro vector fields $\vl_{n}$. A recursive
description for these vector fields were given in \cite{L2} by using
an operator $R$. For any vector field
$\vw = \sum_{n, \alpha} f_{n, \alpha} \tau_{n}(\gamma_{\alpha})$,
define $C(W) := \sum_{n, \alpha, \beta}
         f_{n, \alpha} {\cal C}_{\alpha}^{\beta}
                \tau_{n}(\gamma_{\beta})$, where ${\cal C}$ is the matrix
of multiplication by the first Chern class $c_{1}(M)$
in the ordinary cohomology ring. Define
\[ R(\vw) := \vg * T(\vw) + C(\vw). \]
Then the Virasoro vector fields are given by
\begin{equation} \label{eqn:VirVect}
 \vl_{n} := - R^{n+1}(\vs)
\end{equation}
for $n \geq -1$. One of the nice properties of $\vl_{n}$ is
\[ \bvl_{n} = - \bvx^{n+1} \]
for $n \geq -1$ (cf. \cite[Lemma 4.1]{L2}).
Here, $\bvx^{0}$ is understood as $\bvs$.

\begin{thm} \label{thm:tau-vlf}
For $k \geq -1$,
\[ \overline{\tau_{-}( \vl_{k})} = \sum_{i=1}^{N} u_{i}^{k+1} \vf_{i} -
    \frac{3}{2} (k+1) \bvx^{k}. \]
\end{thm}
{\bf Proof}:
We prove this theorem by induction on $k$.
First note that for $k=-1$, this theorem is precisely \eqref{eqn:tau-S}.
Secondly, by \eqref{eqn:g*ei}, we have
\[
\vg * \bvx^{k}  =  \sum_{i=1}^{N} u_{i}^{k} \vg * \ve_{i}
 =  \frac{1}{2} \, \bvx^{k}  - \sum_{i=1}^{N} u_{i}^{k+1} \vf_{i}
    + \bvx \qp \sum_{i=1}^{N} u_{i}^{k} \vf_{i}.
\]
By \cite[Lemma 4.1 and Theorem 4.8]{L2},
\begin{eqnarray*}
\overline{\tau_{-}(\vl_{k})} &=& \bvx \qp \overline{\tau_{-}(\vl_{k-1})}
        - \bvx^{k} - \vg * \bvx^{k}.
\end{eqnarray*}
Applying the induction hypothesis, we obtain
the desired formula.
$\Box$

\begin{cor} \label{cor:Lkef} For $k \geq -1$,
\[ \vl_{k} =\sum_{i=1}^{N} u_{i}^{k+1} (T(\vf_{i}) - \ve_{i}) -
    \frac{3}{2} (k+1) T\left(\bvx^{k} \right)
    + T^{2}(\tau_{-}^{2}(\vl_{k})). \]
\end{cor}
{\bf Proof}:
This follows from the following special case of \eqref{eqn:WTW}:
\begin{equation} \label{eqn:LkT2}
 \vl_{k} = \bvl_{k} + T(\overline{\tau_{-}(\vl_{k})})
            + T^{2}(\tau_{-}^{2}(\vl_{k}))
\end{equation}
and the fact
\[ \bvl_{k} = - \bvx^{k+1} = -\sum_{i=1}^{N} u_{i}^{k+1} \ve_{i}. \]
$\Box$

\begin{lem} \label{lem:TbS+Hg2}
\[ \gwiitwo{T(\bvs)} = \frac{2}{3} A_{1}(\tau_{-}^{2}(\vl_{0}))
        - \frac{1}{3} \sum_{i=1}^{N} u_{i} B(\ve_{i}, \ve_{i}, \ve_{i}).\]
\end{lem}
{\bf Proof}:
In case
$\vw_{1} = \vw_{2} = \vw_{3} = \ve_{i}$, \eqref{eqn:BP} has a much
simpler form
\begin{equation}
 \gwiitwo{\ve_{i}} -  \gwiitwo{T(\vf_{i})}
= \frac{1}{2} B(\ve_{i}, \ve_{i}, \ve_{i}). \label{eqn:BPei}
\end{equation}
On the other hand, the genus-2 $L_{0}$-constraint has the form
$\gwiitwo{\vl_{0}} = 0$.  Corollary~\ref{cor:Lkef}
and \eqref{eqn:TRR1} then imply the following
\[ 0 = \gwiitwo{\vl_{0}}
    =\sum_{i=1}^{N} u_{i}(\gwiitwo{T(\vf_{i})} - \gwiitwo{\ve_{i}})
        - \frac{3}{2} \gwiitwo{T\left(\bvs \right)}
    + A_{1}(\tau_{-}^{2}(\vl_{0})).
    \]
The lemma then follows from \eqref{eqn:BPei}.
$\Box$

{\bf Proof of Theorem~\ref{thm:F2}}:
Recall that the {\it dilaton vector field} has the
form $\vd = T(\vs)$ (cf. the remark after
\cite[lemma 1.4]{L2}). The standard decomposition of $\vs$ then gives
\[ \vd = T(\bvs + T(\tau_{-}(\vs))) =  T(\bvs) + T^{2}(\tau_{-}(\vs)). \]
By the genus-2 dilaton equation
\[ 2 F_{2} = \gwiitwo{\vd} = \gwiitwo{T(\bvs)}
    + \gwiitwo{T^{2}(\tau_{-}(\vs))}. \]
The theorem then follows from Lemma~\ref{lem:TbS+Hg2} and \eqref{eqn:TRR1}.
$\Box$

{\bf Remark}: A formula for $F_{2}$ under a somewhat weaker
condition was given in \cite[Theorem 5.17]{L2}. The formula given
here is much simpler and much easier to work with than the
corresponding formula in \cite{L2}. Moreover, the proof given here
is also much simpler than the proof in \cite{L2}.

\begin{thm} \label{thm:Lkg2}
If the quantum cohomology is semisimple, then for $k \geq -1$,
\begin{eqnarray*}
 \gwiitwo{\vl_{k}}
     &=&  - \frac{1}{2} \sum_{i=1}^{N} u_{i}^{k+1}
            B(\ve_{i}, \ve_{i}, \ve_{i})
            + \frac{k+1}{4} T(\bvx^{k})
            \left\{ \sum_{i=1}^{N} u_{i}
            B(\ve_{i}, \ve_{i}, \ve_{i}) \right\} \\
    &&  + A_{1}(\tau_{-}^{2}(\vl_{k}))
        - \frac{k+1}{2} \, T(\bvx^{k})
        A_{1} \left(\tau_{-}^{2}(\vl_{0}) \right)
    - \frac{3(k+1)}{4} \, T(\bvx^{k})
            A_{1} \left(\tau_{-}(\vs) \right)
 \end{eqnarray*}
\end{thm}
{\bf Proof}:
By Corollary~\ref{cor:Lkef},
\begin{eqnarray*}
 \gwiitwo{\vl_{k}}
     &=&  - \sum_{i=1}^{N} u_{i}^{k+1}
                \left( \gwiitwo{\ve_{i}} -  \gwiitwo{T(\vf_{i})} \right)
           - \frac{3(k+1)}{2} \, T(\bvx^{k}) F_{2}
       + \gwiitwo{ T^{2} \left( \tau_{-}^{2}(\vl_{k}) \right) }.
 \end{eqnarray*}
The theorem then follows from
applying \eqref{eqn:BPei} to the first term, Theorem~\ref{thm:F2} to
the second term, and \eqref{eqn:TRR1} to the third term.
$\Box$

So far we have used only special cases of \eqref{eqn:BP}. In fact
we can not get more information on the genus-2 generating function from other
cases of this equation. But we can still get some interesting properties of
the complicated tensor $B$ from studying the more general cases of this
equation. We have the following
\begin{lem} \label{lem:PropB}
For $i \neq j \neq k$,
\begin{eqnarray*}
(a) & B(\ve_{i}, \ve_{j}, \ve_{k}) & = 0 \\
(b) & B(\ve_{i}, \ve_{j}, \ve_{j}) & = - B(\ve_{i}, \ve_{i}, \ve_{j}) \\
(c) & B(\ve_{i}, \ve_{i}, \ve_{j}) & =
    \frac{1}{2} T(\ve_{i}) B(\ve_{j}, \ve_{j}, \ve_{j})
    - \frac{1}{2} T(\ve_{j}) B(\ve_{i}, \ve_{i}, \ve_{i})
    + A_{2}(\vf_{j}, \ve_{i}) - A_{2}(\vf_{i}, \ve_{j})
\end{eqnarray*}
\end{lem}
In this lemma, $A_{2}$ is a symmetric 2-tensor which only depends on
genus-0 and genus-1 data. It comes from
a genus-2 equation due to Getzler which takes the following
form (cf. \cite{Ge1}): For any vector fields $\vw_{i}$,
\begin{equation} \label{eqn:TRR2}
 \gwiitwo{T(\vw_{1}) \, T(\vw_{2})} = A_{2}(\vw_{1}, \vw_{2}).
\end{equation}
It was proved in \cite{L2} that this equation follows
from \eqref{eqn:TRR1} and \eqref{eqn:BP}.

{\bf Proof of Lemma~\ref{lem:PropB}}:
It was proved in \cite{L4} that the genus-0 4-point functions satisfy
the following properties: For any $i \neq j \neq k$ and any vector field
$\vw$,
\begin{eqnarray}
&(i)& \gwii{\ve_{i} \, \ve_{j} \, \ve_{k} \, \vw}  = 0, \nonumber \\
&(ii)&  \gwii{ \vw \, \ve_{i} \, \ve_{i} \, \ve_{j}}
        = - \gwii{ \vw \, \ve_{j} \, \ve_{j} \, \ve_{i}},
    \label{eqn:g04ptidem} \\
&(iii)& \gwii{\ve_{i} \, \ve_{i} \, \ve_{j} \, \gua} \ga
    = \vf_{j} \qp \ve_{i} - \vf_{i} \qp \ve_{j}.  \nonumber
\end{eqnarray}
Applying \eqref{eqn:BP} for $\vw_{1}=\ve_{i}$, $\vw_{2}=\ve_{j}$,
$\vw_{2}=\ve_{k}$, we have
\[ B(\ve_{i}, \ve_{j}, \ve_{k})
    = - 2 \gwii{\ve_{i} \, \ve_{j} \, \ve_{k} \, \gua} \gwiitwo{T(\ga)}. \]
Therefore (a) follows from \eqref{eqn:g04ptidem} (i).

Applying \eqref{eqn:BP} for $\vw_{1}=\vw_{2}=\ve_{i}$,
$\vw_{2}=\ve_{j}$, we have
\begin{equation} \label{eqn:BPeiej}
 B(\ve_{i}, \ve_{i}, \ve_{j})
    = - 2 \gwii{\ve_{i} \, \ve_{i} \, \ve_{j} \, \gua} \gwiitwo{T(\ga)}
        - \gwiitwo{\ve_{i} \, T(\ve_{j})}
        + \gwiitwo{\ve_{j} \, T(\ve_{i})}.
\end{equation}
By \eqref{eqn:g04ptidem} (ii), if we interchange $i$ and $j$,
the right hand side of this equation is only changed by a minus sign.
This proves (b).

Since $\nabla_{T(\ve_{j})} \ve_{i} = 0$
for $i \neq j$,
\[ \nabla_{T(\ve_{j})} T(\vf_{i}) = T(\nabla_{T(\ve_{j})} \, \vf_{i})
    = T\left(\gwii{T(\ve_{j}) \, \ve_{i} \, \ve_{i} \, \ve_{i} \, \gua} \ga
        \right).
\]
By \cite[Equation (9)]{L2}, we have
\[ \nabla_{T(\ve_{j})} T(\vf_{i})
     = T\left(\gwii{\ve_{j} \, \ve_{i} \, \ve_{i} \, \gua} \ga \qp \ve_{i}
        \right)
     = T(\vf_{j} \qp \ve_{i}) \]
where the last equality follows from \eqref{eqn:g04ptidem} (iii).
Hence taking derivative of \eqref{eqn:BPei} along
$T(\ve_{j})$, we obtain
\[  \gwiitwo{\ve_{i} \, T(\ve_{j})}
    -  \gwiitwo{T(\vf_{j} \qp \ve_{i})}
    -  \gwiitwo{T(\ve_{j}) \, T(\vf_{i})}
    = \frac{1}{2} T(\ve_{j}) B(\ve_{i}, \ve_{i}, \ve_{i}). \]
Applying \eqref{eqn:TRR2} to the last term on the left hand side of
the above equation, we obtain
\[
  \gwiitwo{\ve_{i} \, T(\ve_{j})}
=   \gwiitwo{T(\vf_{j} \qp \ve_{i})}
    + 3 \gwiitwo{T(\ve_{j} \qp \vf_{i})}
 +  A_{2}(\vf_{i}, \ve_{j})
 + \frac{1}{2} T(\ve_{j}) B(\ve_{i}, \ve_{i}, \ve_{i}).
\]
Plugging this formula into \eqref{eqn:BPeiej}, we obtain
\begin{eqnarray*}
 B(\ve_{i}, \ve_{i}, \ve_{j})
   & = & - 2 \gwii{\ve_{i} \, \ve_{i} \, \ve_{j} \, \gua} \gwiitwo{T(\ga)}
        + 2 \gwiitwo{T(\vf_{j} \qp \ve_{i})}
        - 2 \gwiitwo{T(\vf_{i} \qp \ve_{j})}  \\
   & & +  A_{2}(\vf_{j}, \ve_{i})
    + \frac{1}{2} T(\ve_{i}) B(\ve_{j}, \ve_{j}, \ve_{j})
        -  A_{2}(\vf_{i}, \ve_{j})
    - \frac{1}{2} T(\ve_{j}) B(\ve_{i}, \ve_{i}, \ve_{i})   .
\end{eqnarray*}
The first three terms on the right hand side of this equation are cancelled
with each other due to \eqref{eqn:g04ptidem} (iii). This proves (c).
$\Box$

{\bf Remark}: The reason that the general form of \eqref{eqn:BP}
gives no more information on the genus-2 generating function than
what we can get from \eqref{eqn:BPei} lies behind the proof of this lemma.

\begin{cor} \label{cor:uiBei}
For any integers $m \geq 0$ and $k \geq 0$,
\[
 2 \sum_{i=1}^{N} u_{i}^{m+k} \, B(\ve_{i} , \ve_{i},  \ve_{i})
=  \sum_{i=1}^{N} u_{i}^{m} \, B(\ve_{i}, \ve_{i}, \bvx^{k})
    + \sum_{i=1}^{N} u_{i}^{k} \, B(\ve_{i}, \ve_{i}, \bvx^{m}).
\]
\end{cor}
{\bf Proof}:
Since
 $T(\vw) \, u_{i} = 0$ for any vector field $\vw$,
multiplying Lemma~\ref{lem:PropB} (c) by $u_{i}^{m}u_{j}^{k}$
and summing over $i$ and $j$, we obtain
\begin{eqnarray*}
&& \sum_{i=1}^{N} u_{i}^{m} \, B(\ve_{i}, \ve_{i}, \bvx^{k})
-  \sum_{i=1}^{N} u_{i}^{m+k} \, B(\ve_{i} , \ve_{i},  \ve_{i}) \\
&=&     \frac{1}{2} T(\bvx^{m})
    \sum_{i=1}^{N} u_{i}^{k} \, B(\ve_{i}, \ve_{i}, \ve_{i})
     -  \frac{1}{2} T(\bvx^{k})
    \sum_{i=1}^{N} u_{i}^{m} \, B(\ve_{i}, \ve_{i}, \ve_{i}) \\
&&   + A_{2}( \sum_{i=1}^{N}u_{i}^{k} \vf_{i}, \bvx^{m})
     - A_{2}( \sum_{i=1}^{N}u_{i}^{m} \vf_{i}, \bvx^{k}).
\end{eqnarray*}
Observe that the right hand side of this equation is anti-symmetric
with respect to $m$ and $k$. So the desired formula is obtained by
symmetrizing this equation with respect to $m$ and $k$.
$\Box$

\section{Genus-2 Virasoro conjecture for manifolds with \\ semisimple
quantum cohomology}

By \cite[Theorem 5.9]{L2}, the genus-2 Virasoro conjecture
for any smooth projective variety can be reduced
to the genus-2 $L_{1}$-constraint which have the following form:
\begin{equation} \label{eqn:Virg2L1}
 \gwiitwo{\vl_{1}}= - \frac{1}{2} \sum_{\alpha} b_{\alpha} (1-b_{\alpha})
        \{ \gwiione{\gua \, \ga} + \gwiione{\gua} \gwiione{\ga} \}.
\end{equation}
So to prove Theorem~\ref{thm:g2VirSS}, it suffices to compute
$\gwiitwo{\vl_{1}} = \vl_{1} F_{2}$ and check whether it
coincides with the right hand side of \eqref{eqn:Virg2L1}. There
are two approaches to this problem: The first approach is to take
the formula for $F_{2}$ in Theorem~\ref{thm:F2} and then take the
derivative along $\vl_{1}$. The second approach is to directly use
the formula for $\gwiitwo{\vl_{1}}$ given in
Theorem~\ref{thm:Lkg2}. Since the intermediate results of the
first approach provide more understanding for the genus-2
generating function and the Virasoro vector fields,
we prove Theorem~\ref{thm:g2VirSS} using the
first approach in this section. The second approach will be given in
the appendix.

\subsection{Express everything in terms of idempotents}

\Eqref{eqn:Virg2L1} is given in the flat frame $\{
\gamma_{\alpha} \mid \alpha = 1, \ldots, N\}$, so is the tensor
$A_{1}$ defined after \eqref{eqn:TRR1}. But idempotents has
appeared in the expression of $F_{2}$ as given in Theorem
\ref{thm:F2}. To compare $\vl_{1} F_{2}$ with the formula in
\eqref{eqn:Virg2L1}, we need to re-write both of them using
idempotents only. For this purpose, it is convenient to introduce
the following notation:
\[ z_{i_{1}, \cdots, i_{k}} := \gwii{\ve_{i_{1}} \, \cdots \, \ve_{i_{k}}} \]
and
\[ \phi_{i_{1}, \cdots, i_{k}}
    := \gwiione{\ve_{i_{1}} \, \cdots \, \ve_{i_{k}}}.\]
We will use the following simple fact which was explained in
\cite{L4}: For any tensor $Q$,
\begin{equation} \label{eqn:gagua}
 Q(\ga, \gua, \cdots)
    = \sum_{i=1}^{N} \frac{1}{g_{i}} \, Q(\ve_{i}, \ve_{i}, \cdots).
 \end{equation}
In particular
\begin{equation} \label{eqn:Delta}
\Delta :=  \ga \qp \gua =  \sum_{i=1}^{N} \frac{1}{g_{i}} \,
\ve_{i}.
 \end{equation}
So the prediction of genus-2 $L_{1}$-constraint, i.e.
\eqref{eqn:Virg2L1}, can be re-written as
\begin{equation} \label{eqn:Virg2L1-2}
\gwiitwo{\vl_{1}}
   = - \frac{1}{2} \sum_{i=1}^{N} \frac{1}{g_{i}}
    \left\{ \gwiione{(\vg* \ve_{i}) \, (\vg* \ve_{i})}
        + \gwiione{(\vg* \ve_{i})}^{2} \right\} .
\end{equation}
Using the definition of tensor $B$ in \cite{L2}, we can write down
the precise formula for the function $B(\ve_{i}, \ve_{i},
\ve_{i})$:
\begin{eqnarray}
B(\ve_{i}, \ve_{i}, \ve_{i})
 &=& \sum_{j, k} \frac{1}{g_{j} g_{k}}
    \left\{ \frac{1}{5} z_{iiijk} \phi_{j} \phi_{k}
        - \frac{6}{5} z_{iiij} \phi_{jk} \phi_{k}
        - \frac{6}{5} z_{iijk} \phi_{j} \phi_{ik} \right\} \nonumber \\
&& + \sum_{j} \frac{1}{g_{j}} \left\{ \frac{9}{5}(1-2\delta_{ij})
        \phi_{ij}^{2} -  \frac{6}{5}(1-2\delta_{ij})
        \phi_{iij} \phi_{j} \right\} \nonumber \\
&& + \sum_{j, k} \frac{1}{g_{j} g_{k}}
    \left\{ \frac{1}{120} z_{iiijjk} \phi_{k}
        + \frac{1}{10} z_{iiijk} \phi_{jk}
        - \frac{1}{20} z_{iiij} \phi_{jkk}  \right. \nonumber \\
&& \hspace{60pt} \left. - \frac{3}{40} z_{iijjk} \phi_{ik}
        + \frac{3}{40} z_{ijjk} \phi_{iik}
        - \frac{3}{10} z_{iijk} \phi_{ijk} \right\} \nonumber \\
&& - \sum_{j} \frac{1}{g_{j}} \left\{ \frac{1}{120}
        \phi_{iiij} +  \frac{1}{20}(1-2\delta_{ij})
        \phi_{iijj}  \right\}. \label{eqn:Bi}
\end{eqnarray}

To express $A_{1}$ in terms of idempotents, we will use
decomposition (\ref{eqn:WTW}) and derivatives of genus-0 and
genus-1 topological recursion relations. Recall that the genus-0
topological recursion relation has the following form
\[ \gwii{T(\vw) \, \vv_{1} \, \vv_{2}} = 0 \]
for any vector fields $\vw$ and $\vv_{i}$. Repeatedly taking
derivatives of this relation, we have
\begin{eqnarray}
&&  \gwii{T(\vw) \, \vv_{1} \, \cdots \, \vv_{k+2}} \nonumber  \\
&= &  \sum_{m=1}^{k} \, \, \, \sum_{1 \leq i_{1} < \cdots < i_{m}
\leq k}
    \gwii{\vw \, \vv_{i_{1}} \, \cdots \, \vv_{i_{m}} \, \gua}
        \nonumber \\
&& \hspace{100pt} \cdot
    \gwii{\ga \, \vv_{1} \, \cdots \widehat{\vv_{i_{1}}} \,
        \cdots \widehat{\vv_{i_{2}}} \, \cdots \, \cdots\, \cdots
        \widehat{\vv_{i_{m}}} \, \cdots
        \vv_{k+2}}. \label{eqn:TRRg0}
\end{eqnarray}
The genus-1 topological recursion relation has the form
\[ \gwiione{T(\vw)} = \frac{1}{24} \gwii{\vw \, \ga \, \gua} \]
for any vector field $\vw$. Repeatedly taking derivatives of this
relation, we obtain
\begin{eqnarray}
 && \gwiione{T(\vw) \, \vv_{1} \, \cdots \, \vv_{k}} \nonumber \\
&= &  \sum_{m=1}^{k} \sum_{1 \leq i_{1} < \cdots < i_{m} \leq k}
    \gwii{\vw \, \vv_{i_{1}} \, \cdots \, \vv_{i_{m}} \, \gua}
        \nonumber \\
&& \hspace{100pt} \cdot
    \gwiione{\ga \, \vv_{1} \, \cdots \widehat{\vv_{i_{1}}} \,
        \cdots \widehat{\vv_{i_{2}}} \, \cdots \, \cdots\, \cdots
        \widehat{\vv_{i_{m}}} \, \cdots
        \vv_{k}} \nonumber \\
&& + \frac{1}{24} \gwii{\vw \, \vv_{1}
    \, \cdots \, \vv_{k} \, \gua \, \ga} \label{eqn:TRRg1}
\end{eqnarray}
for all vector fields $\vw$ and $\vv_{i}$,  and all integer $k
\geq 0$. For any vector field $\vw$,  first decomposing it as
\[ \vw = \bvw + T(\overline{\tau_{-}(\vw)}) + T^{2}(\tau_{-}^{2}(\vw)), \]
then using \eqref{eqn:TRRg0} and (\ref{eqn:TRRg1}) to get rid of
the operator $T$, we obtain
\begin{eqnarray*}
 A_{1}(\vw) &=& \frac{7}{10} \gwiione{\gua} \gwiione{(\ga \qp \bvw)} \\
  && + \frac{1}{10} \gwiione{\gua \, (\ga \qp \bvw)}
        -  \frac{1}{240} \gwiione{\bvw \, \Delta}  \\
    &&  + \frac{1}{20}
        \gwiione{\left\{\overline{\tau_{-}(\vw)} \qp \Delta \right\}}
        + \frac{13}{240} \gwii{\bvw \, \ga \, \gua \, \gub}
            \gwiione{\gb}  \\
  && + \frac{1}{1152}
        \gwii{ \vs \, \vs \, \left\{ \overline{\tau_{-}^{2}(\vw)}
            \qp \Delta^{2} \right\} }
        +\frac{1}{1152}
         \gwii{\overline{\tau_{-}(\vw)} \, \Delta \, \gua \, \ga}  \\
    &&  + \frac{1}{480}
        \gwii{\left\{\overline{\tau_{-}(\vw)} \qp \ga \right\}
            \, \gua \, \gb \, \gub}
        + \frac{1}{960}
        \gwii{\bvw \, \ga \, \gua \, \gb \, \gub}.
\end{eqnarray*}
In semisimple case, we can use idempotents to express this tensor
as
\begin{eqnarray}
A_{1}(\vw) &=& \sum_{i} \frac{\left< \vw, \, \ve_{i}
\right>}{g_{i}}
         \left\{ \frac{1}{g_{i}} \left( \frac{7}{10} \phi_{i}^{2}
                + \frac{1}{10} \phi_{ii} \right)
        - \sum_{j}  \frac{1}{240} \frac{1}{g_{j}} \phi_{ij}
                \right. \nonumber  \\
&& \hspace{65pt} \left.     + \sum_{j, k} \left(  \frac{13}{240}
               \frac{1}{g_{j} g_{k}} z_{ijjk} \phi_{k}
               +  \frac{1}{960}
               \frac{1}{g_{j} g_{k}} z_{ijjkk} \right) \right\} \nonumber \\
&&    + \sum_{i}  \frac{\left< \tau_{-}(\vw), \, \ve_{i}
\right>}{g_{i}}
           \left\{  \frac{1}{20} \frac{1}{g_{i}}  \phi_{i}
        + \sum_{j} \frac{1}{480} \frac{1}{g_{i} g_{j}} z_{iijj}
        + \sum_{j, k} \frac{1}{1152}
               \frac{1}{g_{j} g_{k}} z_{ijkk}
             \right\} \nonumber \\
&&    + \sum_{i}  \frac{\left< \tau_{-}^{2}(\vw), \, \ve_{i}
\right>}{g_{i}}
             \,\, \cdot \,\,  \frac{1}{1152}
               \frac{1}{g_{i}}. \label{eqn:A1idem}
\end{eqnarray}

\subsection{Expressing everything by rotation coefficients}
\label{sec:RotCoeff}

Relations among functions $ z_{i_{1}, \cdots, i_{k}}$ and
$\phi_{i_{1}, \cdots, i_{k}}$ are very complicated. It's much
easier to see the relations by introducing rotation coefficients.
On the small phase space, rotation coefficients was introduced by
Dubrovin \cite{D} to study semisimple Frobenius manifolds. A
similar definition for {\it rotation coefficients} on the big
phase space is the following:
\[ r_{ij} := \frac{\ve_{j}}{\sqrt{g_{j}}} \sqrt{g_{i}}. \]
We will briefly review basic properties of rotation coefficients
when they are needed. The readers are referred to \cite{L4} for
more details. First, $(r_{ij})$ is a symmetric matrix. Using
these functions, the operator $\vg *$ is given by
\begin{equation}
 \vg * \ve_{i} = \frac{1}{2} \ve_{i}
+ \sum_{j} (u_{i} - u_{j}) r_{ij} \, \sqrt{\frac{g_{i}}{g_{j}}}
    \, \ve_{j} \label{eqn:G*rot}
\end{equation}
for all $i$. Therefore, the prediction of the genus-2
$L_{1}$-constraint, i.e. \eqref{eqn:Virg2L1-2}, is given by
\begin{eqnarray*}
\gwiitwo{\vl_{1}} &=& - \frac{1}{2} \sum_{i=1}^{N} \frac{1}{g_{i}}
    \left\{ \frac{1}{4}(\phi_{ii} + \phi_{i}^{2})
        + \sum_{j} (u_{i}-u_{j}) r_{ij} \sqrt{\frac{g_{i}}{g_{j}}}
        (\phi_{ij} + \phi_{i} \phi_{j}) \right. \nonumber  \\
&& \left. \hspace{60pt}     + \sum_{j, k} (u_{i}-u_{j})
(u_{i}-u_{k})
         r_{ij} r_{ik} \, \frac{g_{i}}{\sqrt{g_{j} g_{k}}} \,
        (\phi_{jk} + \phi_{j} \phi_{k}) \right\}.
\end{eqnarray*}
Note that the second term is anti-symmetric with respect to $i$
and $j$, so equal to 0 when summing over $i$ and $j$. Define
\[ v_{ij} := (u_{j} - u_{i}) r_{ij}. \]
Then $v_{ij} = - v_{ji}$ for all $i$ and $j$.  The prediction of
the genus-2 $L_{1}$-constraint can now be written as
\begin{eqnarray}
\gwiitwo{\vl_{1}} &=& - \frac{1}{2} \sum_{i=1}^{N} \frac{1}{g_{i}}
    \left\{ \frac{1}{4}(\phi_{ii} + \phi_{i}^{2})
        + \sum_{j, k}
         v_{ij} v_{ik} \, \frac{g_{i}}{\sqrt{g_{j} g_{k}}} \,
        (\phi_{jk} + \phi_{j} \phi_{k}) \right\}.  \label{eqn:Virg2L1-3}
\end{eqnarray}

Covariant derivatives of idempotents are given by
\begin{eqnarray}
 \nabla_{\ve_{i}} \ve_{j}
&=& r_{ij} \left( \sqrt{\frac{g_{j}}{g_{i}}} \, \ve_{i}
            + \sqrt{\frac{g_{i}}{g_{j}}} \, \ve_{j} \right)
  - \delta_{ij} \sum_{k=1}^{N}
         r_{ik} \,  \sqrt{\frac{g_{i}}{g_{k}}} \,\, \ve_{k}
         \label{eqn:covidem}
\end{eqnarray}
for any $i$ and $j$.
To compute the derivatives of rotation coefficients, we define
\[ \theta_{ij} := \frac{1}{u_{j} - u_{i}} \left(
        r_{ij} + \sum_{k} r_{ik} v_{jk} \right) \]
for $i \neq j$. These functions satisfy the following property
\begin{equation} \label{eqn:symmtheta}
 \theta_{ij} + \theta_{ji} = - \sum_{k} r_{ik} r_{jk}
\end{equation}
for any $i \neq j$.
First derivatives of rotation coefficients are given by the formula
\begin{equation}
\ve_{k} r_{ij} = r_{ik} r_{jk}
      + \left\{
         \begin{array}{ll}
       0, & {\rm if} \,\,\, i \neq j \neq k, \\ \\
       \theta_{ij}, & {\rm if} \,\,\, k=i \neq j, \\ \\
       \sqrt{\frac{g_{k}}{g_{i}}} \theta_{ik},
                & {\rm if} \,\,\, i= j \neq k, \\ \\
       - 2 \sum_{l} r_{il}^{2} +
    \sum_{p \neq i} \sqrt{\frac{g_{p}}{g_{i}}} \theta_{pi}
    + \frac{1}{g_{i}} \left<\tau_{-}^{2}(\vs), \, \ve_{i} \right>,
                & {\rm if} \,\,\, i= j = k.
       \end{array}
        \right.
	\label{eqn:ekrij}
\end{equation}
The appearance of $<\tau_{-}^{2}(\vs), \, \ve_{i}>$ in the last
equation is a typical big phase space phenomenon. This term
vanishes on the small phase space, but is in general not zero on
the big phase space.  To compute higher order derivatives of
rotation coefficients,  we also need the following formula:
\begin{equation}
 \ve_{j} \, \, <\tau_{-}^{k}(\vs), \, \ve_{i}> \,\, = \,\,
    \delta_{i j} <\tau_{-}^{k+1}(\vs), \, \ve_{i}>
    + <\tau_{-}^{k}(\vs), \, \nabla_{\ve_{j}} \ve_{i}> 
	\label{eqn:eitauS}
\end{equation}
for all $k \geq 0$. This formula follows from
\eqref{eqn:innermetric}. We can repeatedly apply the above
formulas to compute higher order derivatives of rotation
coefficients in terms of functions $g_{i}$, $u_{i}$, $r_{ij}$, and
$<\tau_{-}^{k}(\vs), \, \ve_{i}>$ with $k \geq 2$. When we
compute $k$-th order derivatives of $r_{ij}$, we might encounter
$<\tau_{-}^{k+1}(\vs), \, \ve_{i}>$. For most purposes (for
example, for the proof of genus-1 and genus-2 Virasoro
conjecture) these terms will not affect the final results.
In this paper, we only need derivatives of rotation coefficients up
to order 3. It will be convenient to introduce the following functions
when computing the second and third order derivatives:
\begin{eqnarray*}
 \Omega_{ij} &:=& \frac{1}{u_{j}-u_{i}} \left\{
          \theta_{ij} - \theta_{ji} + \sum_{k, l}
         r_{il} r_{jk} v_{kl} \right\}, \\
\Lambda_{ij} &:=&  \frac{1}{u_{j}- u_{i}} \left\{
       3 \Omega_{ij} - \sum_{k}(u_{k}-u_{j}) \theta_{ik} \theta_{jk}
       -  \left( r_{ii} + \sum_{k} r_{ik} v_{ik} \right) \theta_{ji}
       \right. \nonumber \\
&& \hspace{80pt} \left.
       -  \left( r_{jj} + \sum_{k} r_{jk} v_{jk} \right) \theta_{ij}  \right\}
\end{eqnarray*}
for $i \neq j$. These functions have the property
\[ \Omega_{ij} = \Omega_{ji} \hspace{30pt} {\rm and} \hspace{30pt}
    \Lambda_{ij} + \Lambda_{ji} = \sum_{k} \theta_{ik} \theta_{jk}
\]
for all $i \neq j$. They arise naturally in the second and third
order derivatives of rotation coefficients because
\[ \ve_{j} \theta_{ij} = \left( r_{jj} - \sqrt{\frac{g_{j}}{g_{i}}} r_{ij}
                       \right) \theta_{ij}
               - \Omega_{ij} \]
and
\[ \ve_{i} \Omega_{ij}
  = \theta_{ji} \left( \sum_{k} r_{ik}^{2}
         + \sqrt{\frac{g_{i}}{g_{j}}}  \theta_{ij}
         + \sum_{k \neq i} \sqrt{\frac{g_{k}}{g_{i}}}  \theta_{ki}
     \right) + \Lambda_{ij} \]
for $i \neq j$. We might consider $\theta_{ij}$, $\Omega_{ij}$, and
$\Lambda_{ij}$ as functions having poles of order 1, 2, and 3 respectively
in terms of $u_{1}, \ldots, u_{n}$.

We can assign $g_{i}$ with degree 0, $u_{i}$ with degree $-1$,
and $k$-th order derivatives of $r_{ij}$ along directions of idempotents
with degree $k+1$. Then most expressions in this paper are
homogeneous of a fixed degree. For example, functions $z_{i_{1},
\cdots, i_{k}}$ have degree $k-3$ for $k \geq 3$, $\phi_{i_{1},
\cdots, i_{k}}$ have degree $k$, 
$\left< \tau_{-}^{k}(\vs), \,  \ve_{i} \right>$
has degree $k$, 
$\left< \tau_{-}^{k}(\vl_{m}), \,  \ve_{i} \right>$
has degree $k-m-1$, 
$B(\ve_{i}, \ve_{i}, \ve_{i})$
has degree 4,  $\gwiitwo{\vl_{1}}$ has degree 2, and 
$F_{2}$ has degree 3. In general, we expect $F_{g}$ to have degree $3(g-1)$
for all $g$.

To express $F_{2}$ in terms of rotation coefficients, we need
first to use rotation coefficients to describe vector fields
$\overline{\tau_{-}^{k}(\vl_{0})}$.
A recursion formula for
$\overline{\tau_{-}^{k}(\vl_{m})}$ was given in \cite[Theorem 4.8]{L2}.
In the semisimple case, using equations (\ref{eqn:Frob}),
(\ref{eqn:G*inner}) and (\ref{eqn:G*rot}), this recursion relation can be
written as
\begin{eqnarray}
\left<\tau_{-}^{k+1}(\vl_{m}), \, \ve_{i} \right>
&=& u_{i} \left<\tau_{-}^{k+1}(\vl_{m-1}), \, \ve_{i} \right>
    + (k + \frac{3}{2}) \left<\tau_{-}^{k}(\vl_{m-1}), \, \ve_{i} \right>
    \nonumber \\
&&    + \sum_{j} v_{ij} \sqrt{\frac{g_{i}}{g_{j}}}
            \left<\tau_{-}^{k}(\vl_{m-1}), \, \ve_{j} \right>.
                \label{eqn:LmRec}
\end{eqnarray}
In particular, since $\vl_{-1} = - \vs$, we have
\begin{eqnarray}
< \tau_{-}^{k} (\vl_{0}), \, \ve_{i}> &=&
    - u_{i} < \tau_{-}^{k}(\vs), \, \ve_{i}>
    - \frac{2k+1}{2} <\tau_{-}^{k-1}(\vs), \, \ve_{i}>   \nonumber \\
&&  - \sum_{j} v_{ij}  \sqrt{\frac{g_{i}}{g_{j}}}
    <\tau_{-}^{k-1}(\vs), \, \ve_{j}>  \label{eqn:L0idem}
\end{eqnarray}
for $k \geq 1$. We will keep $ <\tau_{-}^{k}(\vs), \, \ve_{i}>$
with $k \geq 2$ in our computations. But for $k=0$ and $k=1$, we
will use (cf. \cite{L4}),
\begin{equation}  < \vs, \, \ve_{i}> = g_{i} \hspace{20pt} {\rm and}  \hspace{20pt}
   < \tau_{-}(\vs), \, \ve_{i}> = \sum_{j} r_{ij} \sqrt{g_{i} g_{j}} \, \, .
   \label{eqn:Sidem}
\end{equation}

We now describe how to represent $z_{i_{1}, \ldots, i_{k}}$ and
$\phi_{i_{1}, \ldots, i_{k}}$ in terms of rotation coefficients.
In \cite{L4}, it was proved that genus-0 4-point functions have
the following property: For $i \neq j$,
\begin{eqnarray} &(i)& z_{iiii}
    = - g_{i} r_{ii}, \nonumber \\
&(ii)& z_{jiii}
    = - z_{jjii}
    = - \sqrt{g_{i} g_{j}} \,  r_{ij}, \nonumber \\
&(iii)& z_{ijkl} = 0 \hspace{20pt} {\rm otherwise}. \label{eqn:z4rot}
\end{eqnarray}
It was also proved that genus-1 1-point functions are given by
\begin{eqnarray}
24 \phi_{i} &=& - 12 \sum_{j}  r_{ij} v_{ij}
    - \sum_{j}  \sqrt{\frac{g_{i}}{g_{j}}} \, r_{ij}. \label{eqn:phi1rot}
\end{eqnarray}
Note that for any positive integer $k$,
\begin{eqnarray}
 \gwiig{\ve_{i_{1}} \ve_{i_{2}} \cdots \ve_{i_{k+1}}}
    &=& \ve_{i_{k+1}} \gwiig{\ve_{i_{1}} \cdots \ve_{i_{k}}}
        - \sum_{j=1}^{k} \gwiig{\ve_{i_{1}} \cdots
            \left(\nabla_{\ve_{i_{k+1}}} \ve_{i_{j}}\right)
        \cdots \ve_{i_{k}}} \nonumber \\
    &=&  \ve_{i_{k+1}} \gwiig{\ve_{i_{1}} \cdots \ve_{i_{k}}}
        - \left( \sum_{j=1}^{k} r_{i_{j}, i_{k+1}}
    \sqrt{\frac{g_{i_{k+1}}}{g_{i_{j}}}}
            \right)\gwiig{\ve_{i_{1}} \cdots \ve_{i_{k}}} \nonumber \\
    & &  -  \sum_{j=1}^{k} r_{i_{j}, i_{k+1}}
        \sqrt{\frac{g_{i_{j}}}{g_{i_{k+1}}}}
             \gwiig{\ve_{i_{1}} \cdots \widehat{\ve_{i_{j}}} \cdots \ve_{i_{k}}
             \ve_{i_{k+1}}}   \nonumber \\
    & &    + \sum_{j=1}^{k} \delta_{i_{k+1}, i_{j}}
            \sum_{p} r_{p, i_{k+1}} \sqrt{\frac{g_{i_{K+1}}}{g_{p}}}
            \gwiig{\ve_{i_{1}} \cdots \widehat{\ve_{i_{j}}}
          \cdots \ve_{i_{k}} \ve_{p}}. \label{eqn:CorrRecur}
\end{eqnarray}
Repeatedly using this formula, we can express all functions
$z_{i_{1} \cdots i_{k}}$, with $k \geq 4$, and $\phi_{i_{1}
\cdots i_{k}}$, with $k \geq 1$, in terms of rotation
coefficients. For example, genus-1 2-point functions are given by
\begin{eqnarray}
24 \phi_{ii} &=&   12 r_{ii}^{2}
     - \frac{1}{g_{i}} \left< \tau_{-}^{2}(\vs), \, \ve_{i} \right>
     + \sum_{j} \left\{ r_{ij}^{2}
         \left( -10 + \frac{g_{i}}{g_{j}} \right)
           + 24 r_{ii} r_{ij} v_{ij} \right\} \nonumber \\
&&    - \sum_{j,k} \left( 12 r_{ij} r_{jk} v_{jk}
\sqrt{\frac{g_{i}}{g_{j}}}
     + r_{ij} r_{jk} \sqrt{\frac{g_{i}}{g_{k}}} \right)
    - \sum_{j \neq i} \left( \theta_{ij}  \sqrt{\frac{g_{i}}{g_{j}}}
    +  \theta_{ji}  \sqrt{\frac{g_{j}}{g_{i}}} \right) \label{eqn:phi2iirot}
\end{eqnarray}
for all $i$, and
\begin{eqnarray}
24 \phi_{ij} &=& 12 r_{ij}^{2} + \sum_{k} \left\{ r_{ik} r_{jk}
\frac{\sqrt{g_{i} g_{j}}}{g_{k}}
           + 12  r_{ij} r_{ik} v_{ik}  \sqrt{\frac{g_{j}}{g_{i}}}
       + 12 r_{ij} r_{jk} v_{jk}  \sqrt{\frac{g_{i}}{g_{j}}}
       \right\} \nonumber \\
&&        - \left( \theta_{ij}  \sqrt{\frac{g_{j}}{g_{i}}}
                    +  \theta_{ji}  \sqrt{\frac{g_{i}}{g_{j}}} \right)
                    \label{eqn:phi2ijrot}
\end{eqnarray}
for $i \neq j$. When $k$ becomes larger, the formula for
$z_{i_{1} \cdots i_{k}}$ and $\phi_{i_{1} \cdots i_{k}}$ becomes
more complicated. It is not illuminating to write them out here.
But one should notice that to express $F_{2}$ in terms of rotation
coefficients, we only need $z_{i_{1} \cdots i_{k}}$ for $4 \leq k
\leq 6$ and $\phi_{i_{1} \cdots i_{k}}$ for $1 \leq k \leq 4$, which can
be obtained by taking derivaties of equations (\ref{eqn:z4rot}),
 (\ref{eqn:phi2iirot}), (\ref{eqn:phi2ijrot}) twice.
Combining with the formula in Theorem~\ref{thm:F2}, and equations
(\ref{eqn:Bi}),  (\ref{eqn:A1idem}), (\ref{eqn:L0idem}), (\ref{eqn:Sidem}),
a lengthy but straightforward computation shows
the following:
\begin{thm} \label{thm:F2rot}
For any manifold with semisimple quantum cohomology, the genus-2
generating function for the Gromov-Witten invariants is given by the
following formula:
\begin{eqnarray*}
&& 5760 \,\, F_{2} \nonumber \\
&=& - 5 \sum_{i} \frac{1}{g_{i}^{2}}
            \left< \tau_{-}^{3}(\vs), \, \ve_{i} \right> 
    + \sum_{i} \sum_{j \neq i} 5 \Omega_{ij} \left(
    \frac{1}{g_{j}}\sqrt{\frac{g_{i}}{g_{j}}}
    - \frac{1}{\sqrt{g_{i} g_{j}}} \right)  \nonumber \\
&& + \sum_{i} \frac{1}{g_{i}}
            \left< \tau_{-}^{2}(\vs), \, \ve_{i} \right> \left\{
        24 r_{ii} \frac{1}{g_{i}}
        + \sum_{j} \left(
           5 r_{ij} \frac{1}{g_{j}}\sqrt{\frac{g_{i}}{g_{j}}}
           + 144 r_{ij} v_{ij} \frac{1}{g_{i}}
            \right)
        \right\}  \nonumber \\
&&  + \sum_{i} \sum_{j \neq i} \theta_{ij} \left\{
        - 24 r_{ii} \frac{1}{g_{i}}\sqrt{\frac{g_{j}}{g_{i}}}
        + 200 r_{ij} \frac{1}{g_{j}}
    \right. \nonumber \\
&& \hspace{66pt} \left.
        + \sum_{k} \left[
            r_{ik} v_{ik} \left( 120 \frac{1}{\sqrt{g_{i} g_{j}}}
               - 144  \frac{1}{g_{i}}\sqrt{\frac{g_{j}}{g_{i}}}
                \right)
            + r_{jk} v_{ik} \left( 85 \frac{1}{g_{i}}
                + 45 \frac{1}{g_{j}} \right) \right]
        \right\} \nonumber \\
&& - \sum_{i}
     576 r_{ii}^{3} \frac{1}{g_{i}} - \sum_{i} 576 \frac{1}{g_{i}}
     \left( \sum_{j} r_{ij} v_{ij} \right)^{3}
     \nonumber \\
&& + \sum_{i,j} \left\{
     480 r_{ij}^{3} \frac{1}{\sqrt{g_{i} g_{j}}}
    - 23 r_{ii} r_{ij}^{2} \frac{1}{g_{i}}
    - 1728 r_{ii}^{2} r_{ij} v_{ij} \frac{1}{g_{i}}
    \right\} \nonumber \\
&& + \sum_{i,j, k} \left\{
    - 24 r_{ii} r_{ik} r_{jk} \frac{1}{g_{i}}\sqrt{\frac{g_{j}}{g_{i}}}
    + 115 r_{ij} r_{ik} r_{jk} \frac{1}{g_{i}}
      \right.    \nonumber \\
&& \hspace{40pt} \left.
    + 1452 r_{ik}^{2} (r_{ij} v_{ij}) \frac{1}{g_{i}}
    - 1728 r_{ii} \frac{1}{g_{i}}
     \left( r_{ij} v_{ij} \right)
        \left( r_{ik} v_{ik} \right)
         \right\} \nonumber \\
&& + \sum_{i, j, k, l} \left\{
        120 r_{ik} r_{jk} (r_{il} v_{il})
                \frac{1}{\sqrt{g_{i} g_{j}}}
     - 144 r_{ij} r_{il} (r_{jk} v_{jk}) 
	\frac{1}{g_{j}}\sqrt{\frac{g_{l}}{g_{j}}}
      \right.    \nonumber \\
&& \hspace{40pt} \left.
     - 40 r_{ik} r_{jk} r_{il} v_{jl} \frac{1}{g_{i}}
     + 720 r_{ij} (r_{ik} v_{ik})(r_{jl} v_{jl})
               \frac{1}{\sqrt{g_{i} g_{j}}}
     \right\}. 
\end{eqnarray*}
\end{thm}
This formula is much simpler than the formulas for expressing 4-point 
genus-1 functions
$\phi_{iiij}$ and $\phi_{iijj}$ in terms of rotation coefficients
(which were omitted here).
We also note that the right hand side of this formula only depends on
genus-0 data.
In \cite[p157-160]{DZ3}, a three-page formula of $F_{2}$ for semisimple
Frobenius manifolds was derived under the assumption that
$F_{2}$ satisfies the Virasoro constraints.
Comparitively, Theorem~\ref{thm:F2rot} gives a much 
simpler formula for $F_{2}$. 

Note that when computing 4-point genus-1 functions $\phi_{iiij}$ and $\phi_{iijj}$
using equations (\ref{eqn:phi1rot}) and (\ref{eqn:CorrRecur}), we will encounter
third order derivatives of rotation coefficients. This causes that
 $B(\ve_{i}, \ve_{i}, \ve_{i})$ contains third order poles of the form
\[ \frac{1}{576} \sum_{j \neq i} \left\{
         \Lambda_{ji} \left( \frac{1}{\sqrt{g_{i} g_{j}}}
                - \frac{1}{g_{j}} \sqrt{\frac{g_{i}}{g_{j}}} \right)
        -  \Lambda_{ij} \left( \frac{1}{\sqrt{g_{i} g_{j}}}
      - \frac{1}{g_{i}} \sqrt{\frac{g_{j}}{g_{i}}} \right)
                \right\}. \]
When interchanging $i$ and $j$, the second term in this expression
is precisely the first term with an opposite
sign.  Multiplying this expression by $u_{i}$ and summing over $i$, we obtain
\[ \frac{1}{576} \sum_{i} \sum_{j \neq i} (u_{i} - u_{j})
         \Lambda_{ji} \left( \frac{1}{\sqrt{g_{i} g_{j}}}
                - \frac{1}{g_{j}} \sqrt{\frac{g_{i}}{g_{j}}} \right). \]
By the definition of $\Lambda_{ij}$, $(u_{i} - u_{j}) \Lambda_{ji}$  only has second order
poles. Therefore third order poles do not appear in the formula for $F_{2}$ in
Theorem~\ref{thm:F2rot}. Similar observations have also been used to simply terms with
first and second order poles when computing $F_{2}$.

Also note that the 4-point genus-1 functions $\phi_{iiii}$ will produce a term
$ - \frac{1}{576} \left< \tau_{-}^{4}(\vs), \, \ve_{i} \right>$
in $B(\ve_{i}, \ve_{i}, \ve_{i})$. After multiplied by $u_{i}$ and summed over $i$,
this term will be cancelled with the corresponding
term produced by $2 A_{1}(\tau_{-}^{2}(\vl_{0}))$. Therefore the formula for 
$F_{2}$ in Theorem~\ref{thm:F2rot}
only contains $\left< \tau_{-}^{k}(\vs), \, \ve_{i} \right>$ for $k=2, 3$.

\subsection{Action of $\vl_{n}$}

In this section, we discuss the action of Virasoro vector fields
on functions $u_{i}$, $g_{i}$, $r_{ij}$, and $\left< \tau_{-}^{k}(\vs), \, \ve_{i} \right>$.
Although for the proof of Theorem~\ref{thm:g2VirSS} we only need to know the action
of $\vl_{1}$, we will discuss the action for all $\vl_{m}$ since
they may be needed for the study of higher genus Gromov-Witten invariants.

Recall for any vector field $\vw$, we have the standard decomposition
\[ \vw = \bvw + T(\tau_{-}(\vw)). \]
In previous sections, we have given formulas for the action of $\ve_{i}$
on functions $u_{i}$, $g_{i}$, $r_{ij}$, and $\left< \tau_{-}^{k}(\vs), \, \ve_{i} \right>$.
Since for any vector filed $\vw$, its primary projection $\bvw$ has the decomposition
\[ \bvw = \sum_{i} \frac{\left<\vw, \, \ve_{i} \right>}{g_{i}} \ve_{i}, \]
the action of $\bvw$ is thus well understood. In \cite{L4}, we have shown the following
formula
\begin{eqnarray*}
T(\vw) u_{i} &=& 0, \\
T(\vw) g_{i} &=& - 2 \left<\vw, \, \ve_{i} \right>, \\
T(\vw) r_{ij} &=& \delta_{ij} \left\{ - \frac{\left<\tau_{-}(\vw), \, \ve_{i} \right>}{g_{i}}
        +  \sum_{k} r_{ik} \frac{1}{\sqrt{g_{i} g_{k}}} \left<\vw, \, \ve_{k} \right>
            \right\}
\end{eqnarray*}
for any vector field $\vw$. Therefore for any vector field $\vw$, we have the followings
formula
\begin{eqnarray}
\vw u_{i} &=& \frac{1}{g_{i}} \left<\vw, \, \ve_{i} \right>, \nonumber \\
\vw g_{i} &=& - 2 \left<\tau_{-}(\vw), \, \ve_{i} \right>
                + \sum_{j} 2 r_{ij} \sqrt{\frac{g_{i}}{g_{j}}}
                    \left<\vw, \, \ve_{j} \right>, \nonumber \\
\vw r_{ij} &=& \frac{1}{g_{i}} \left<\vw, \, \ve_{i} \right> \theta_{ij}
            + \frac{1}{g_{j}} \left<\vw, \, \ve_{j} \right> \theta_{ji}
            + \sum_{k} \frac{1}{g_{k}} \left<\vw, \, \ve_{k} \right> r_{ik} r_{jk}
            \hspace{20pt} {\rm for} \hspace{6pt} i \neq j, \nonumber \\
\vw r_{ii} &=& \sum_{j \neq i} \left\{
                \frac{1}{g_{j}} \left<\vw, \, \ve_{j} \right> \theta_{ij}
            + \frac{1}{g_{i}} \left<\vw, \, \ve_{i} \right> \theta_{ji}
                \right\} \sqrt{\frac{g_{j}}{g_{i}}} \nonumber \\
&&          - \frac{1}{g_{i}} \left<\tau_{-}^{2}(\vw), \, \ve_{i} \right>
            + \frac{1}{g_{i}^{2}} \left< \vw, \, \ve_{i} \right>
                \left<\tau_{-}^{2} (\vs), \, \ve_{i} \right>
            + \sum_{j} \frac{r_{ij}}{\sqrt{g_{i} g_{j}}}
                \left<\tau_{-}(\vw), \, \ve_{j} \right>
                \nonumber \\
&&          + \sum_{j} r_{ij}^{2} \left\{
                \frac{1}{g_{j}} \left<\vw, \, \ve_{j} \right>
                - 2 \frac{1}{g_{i}} \left<\vw, \, \ve_{i} \right>
                \right\}.
\end{eqnarray}

For the Virasoro vector fields $\vl_{m}$, we have
\[ \left<\vl_{m}, \, \ve_{i} \right> = - u_{i}^{m+1} g_{i} \]
since $\bvl_{m} = - \bvx^{m+1}$.
The recursion formula (\ref{eqn:LmRec}) implies that
\begin{equation} \left<\tau_{-}(\vl_{m}), \, \ve_{i} \right>
    = - \frac{3}{2}(m+1) u_{i}^{m} g_{i}
      - \sum_{j} u_{j}^{m+1} r_{ij} \sqrt{g_{i} g_{j}}
    \label{eqn:tauLmRot}
\end{equation}
and
\begin{eqnarray}
\left<\tau_{-}^{2}(\vl_{m}), \, \ve_{i} \right>
&=& - u_{i}^{m+1} \left<\tau_{-}^{2}(\vs), \, \ve_{i} \right>
    - \frac{1}{8} m (11m + 19) u_{i}^{m-1} g_{i} \nonumber \\
&&    - \sum_{j} r_{ij} \sqrt{g_{i} g_{j}}
            \left( 2m u_{i}^{m} + \frac{1}{2}(3m+7) u_{j}^{m}
                    - \sum_{p=0}^{m} u_{i}^{p} u_{j}^{m-p} \right) \nonumber \\
&&    - \sum_{j,k} v_{ij} r_{jk} \sqrt{g_{i} g_{k}}
            \sum_{p=0}^{m} u_{i}^{p} u_{k}^{m-p}.
            \label{eqn:tau2LmRot}
\end{eqnarray}
These formulas enable us to compute $\vl_{m} u_{i}$, $\vl_{m} g_{i}$, and
$\vl_{m} r_{ij}$. When computing $\vl_{m} r_{ij}$, we will encounter first order poles
of the type $u_{i}^{m+1} \theta_{ij} + u_{j}^{m+1} \theta_{ji}$. By \eqref{eqn:symmtheta},
we have
\[ u_{i}^{m+1} \theta_{ij} + u_{j}^{m} \theta_{ji}
    = (u_{i}^{m+1} - u_{j}^{m+1}) \theta_{ij} - u_{j}^{m+1} \sum_{k} r_{ik} r_{jk}.\]
Factoring out $u_{i} - u_{j}$ from $u_{i}^{m+1} - u_{j}^{m+1}$ and using the definition
of $\theta_{ij}$, we have
\[ u_{i}^{m+1} \theta_{ij} + u_{j}^{m+1} \theta_{ji}
    = - r_{ij} \sum_{p=0}^{m} u_{i}^{p} u_{j}^{m-p}
        - \sum_{k} r_{ik} r_{jk} \left\{ u_{j}^{m+1} +
            (u_{k} - u_{j}) \sum_{p=0}^{m} u_{i}^{p} u_{j}^{m-p} \right\}. \]
Similarly,
\[ u_{j}^{m+1} \theta_{ij} + u_{i}^{m+1} \theta_{ji}
    =  r_{ij} \sum_{p=0}^{m} u_{i}^{p} u_{j}^{m-p}
        + \sum_{k} r_{ik} r_{jk} \left\{ - u_{i}^{m+1} +
            (u_{k} - u_{j}) \sum_{p=0}^{m} u_{i}^{p} u_{j}^{m-p} \right\}. \]
Therefore first order poles $\theta_{ij}$ will not appear in $\vl_{m} r_{ij}$ for
all $i$ and $j$. A straightforward computation shows the following
\begin{lem} \label{lem:actLm}
For $i \neq j$,
\begin{eqnarray*}
\vl_{m} u_{i} &=& - u_{i}^{m+1}, \\
\vl_{m} g_{i} &=& 3 (m+1) u_{i}^{m} g_{i}, \\
\vl_{m} r_{ij} &=& r_{ij} \sum_{p=0}^{m} u_{i}^{p} u_{j}^{m-p}
        + \sum_{k} r_{ik} r_{jk} \left\{ u_{j}^{m+1} - u_{k}^{m+1}
                + (u_{k}-u_{j}) \sum_{p=0}^{m} u_{i}^{p} u_{j}^{m-p}
                \right\}, \\
\vl_{m} r_{ii} &=& \frac{1}{8} m (11m+19) u_{i}^{m-1} + (m+1) u_{i}^{m} r_{ii} \\
&&        + \sum_{k} r_{ik} \sqrt{\frac{g_{k}}{g_{i}}} \left( 2m u_{i}^{m}
                - \sum_{p=1}^{m} 2 u_{i}^{p} u_{k}^{m-p} \right) \\
&&        + \sum_{k} r_{ik}^{2} \left\{ (m+1) u_{i}^{m} u_{k}
                    - m u_{i}^{m+1} - u_{k}^{m+1} \right).
\end{eqnarray*}
\end{lem}
From this lemma, we see that the action of $\vl_{m}$ dereases the degree
by $m$. 

We now consider $\vl_{m} \left< \tau_{-}^{k}(\vs), \, \ve_{i} \right>$. We first note that
for any vector field $\vw$,
\[ \nabla_{\vw} \, \tau_{-}^{k}(\vs)
    = \tau_{-}^{k}(\nabla_{\vw} \vs) = - \tau_{-}^{k+1}(\vw) \]
for any $k \geq 0$. So by \eqref{eqn:innermetric},
\begin{equation}
 \vw \left< \tau_{-}^{k}(\vs), \, \ve_{i} \right>
 = - \left< \tau_{-}^{k+1}(\vw), \, \ve_{i} \right>
    + \left< \tau_{-}^{k+1}(\vs), \, \vw \qp \ve_{i} \right>
    + \left< \tau_{-}^{k}(\vs), \, \nabla_{\vw} \ve_{i} \right>.
    \label{eqn:dertauei}
\end{equation}
For the Virasoro vector field $\vl_{m}$, using the standard decomposition, we have
\[ \nabla_{\vl_{m}} \, \ve_{i}
    = - \nabla_{\bvx^{m+1}} \, \ve_{i} + \nabla_{T(\tau_{-}(\vl_{m}))} \, \ve_{i}
    = - \sum_{j} u_{j}^{m+1} \nabla_{\ve_{j}} \ve_{i}
        - \tau_{-}(\vl_{m}) \qp \ve_{i}. \]
By equations (\ref{eqn:covidem}) and  (\ref{eqn:tauLmRot}), we obtain the
formula
\begin{equation}
\nabla_{\vl_{m}} \, \ve_{i} =
    \frac{3}{2} (m+1) u_{i}^{m} \ve_{i}
    + \sum_{j} (u_{i}^{m+1} - u_{j}^{m+1}) r_{ij} \sqrt{\frac{g_{i}}{g_{j}}} \ve_{j}
\label{eqn:deridemLm}
\end{equation}
So by \eqref{eqn:dertauei}, we have
\begin{eqnarray}
 \vl_{m} \left< \tau_{-}^{k}(\vs), \, \ve_{i} \right>
&=& - \left< \tau_{-}^{k+1}(\vl_{m}), \, \ve_{i} \right>
    - u_{i}^{m+1} \left< \tau_{-}^{k+1}(\vs), \, \ve_{i} \right>
    + \frac{3}{2} (m+1) u_{i}^{m} \left< \tau_{-}^{k}(\vs), \, \ve_{i} \right>
        \nonumber \\
&&    + \sum_{j} (u_{i}^{m+1} - u_{j}^{m+1}) r_{ij} \sqrt{\frac{g_{i}}{g_{j}}}
        \left< \tau_{-}^{k}(\vs), \, \ve_{j} \right>.
    \label{eqn:dertaueiLm}
\end{eqnarray}
Note the first term on the right hand side, i.e.
$\left< \tau_{-}^{k+1}(\vl_{m}), \, \ve_{i} \right>$ can be computed recursively
using \eqref{eqn:LmRec}.

\subsection{Computing $\vl_{1} F_{2}$}

We are now ready to compute $\vl_{1} F_{2}$ where $F_{2}$ is given by
Theorem~\ref{thm:F2rot}.
First observe that for $\vl_{1}$, Lemma~\ref{lem:actLm} has a simpler form
\begin{lem} \label{lem:L1}
For all $i$ and $j$,
\begin{eqnarray*}
&(i)& \vl_{1} u_{i} = - u_{i}^{2}, \\
&(ii)& \vl_{1} g_{i} = 6 u_{i} g_{i}, \\
&(iii)& \vl_{1} r_{ij} = \frac{15}{4} \delta_{ij} + (u_{i}+u_{j}) r_{ij}
        - \sum_{k} v_{ik} v_{jk}.
\end{eqnarray*}
\end{lem}
Using this lemma and the definition of $v_{ij}$,
 $\theta_{ij}$ and $\Omega_{ij}$, we obtain
for $i \neq j$,
\begin{equation}
\vl_{1} v_{ij} = (u_{i} - u_{j}) \sum_{k} v_{ik} v_{jk},
        \label{eqn:L1vij}
\end{equation}
\begin{equation}
\vl_{1} \theta_{ij} = (3 u_{i} + u_{j} ) \theta_{ij} - \frac{11}{4} r_{ij}
            + \sum_{k,l} r_{ik} v_{jl} v_{kl},
            \label{eqn:L1theta}
\end{equation}
and
\begin{equation}
\vl_{1} \Omega_{ij} = 3 ( u_{i} + u_{j} ) \Omega_{ij}
            - \frac{11}{4} \sum_{k} r_{ik} r_{jk}
            + \sum_{k, l, p} r_{il} r_{jk} v_{kp} v_{lp}.
            \label{eqn:L1Omega}
\end{equation}
In particular
\begin{eqnarray}
&& \vl_{1} \left\{ \Omega_{ij} \left( \frac{1}{\sqrt{g_{i} g_{j}}}
        - \frac{1}{g_{j}} \sqrt{\frac{g_{i}}{g_{j}}}  \right) \right\} \nonumber \\
&=& 6 \left( \theta_{ij} - \theta_{ji}
         + \sum_{k, l} r_{il} r_{jk} v_{kl} \right)
                     \frac{1}{g_{j}} \sqrt{\frac{g_{i}}{g_{j}}} \nonumber \\
&&   + \left( - \frac{11}{4} \sum_{k} r_{ik} r_{jk}
            + \sum_{k, l, p} r_{il} r_{jk} v_{kp} v_{lp} \right)
         \left( \frac{1}{\sqrt{g_{i} g_{j}}}
        - \frac{1}{g_{j}} \sqrt{\frac{g_{i}}{g_{j}}}  \right).
        \label{eqn:L12ndpoleF2}
\end{eqnarray}
Therefore $\vl_{1} F_{2}$ does not contain second order poles.
Moreover by \eqref{eqn:LmRec},
\begin{eqnarray}
\left< \tau_{-}^{m}(\vl_{1}), \, \ve_{i} \right> &=&
    - u_{i}^{2} < \tau_{-}^{m}(\vs), \, \ve_{i}>
    - (2m+1) u_{i} <\tau_{-}^{m-1}(\vs), \, \ve_{i}>  \nonumber  \\
&&  - \sum_{j} (u_{j}^{2}- u_{i}^{2}) r_{ij} \sqrt{\frac{g_{i}}{g_{j}}}
    <\tau_{-}^{m-1}(\vs), \, \ve_{j}>   \nonumber  \\
&&  - ( m^{2} - \frac{1}{4}) <\tau_{-}^{m-2}(\vs), \, \ve_{i}>
    - 2m \sum_{j} v_{ij} \sqrt{\frac{g_{i}}{g_{j}}}
    <\tau_{-}^{m-2}(\vs), \, \ve_{j}>   \nonumber  \\
&&  -  \sum_{j,k} v_{ij} v_{jk} \sqrt{\frac{g_{i}}{g_{k}}}
    <\tau_{-}^{m-2}(\vs), \, \ve_{k}>
    \label{eqn:taukL1ei}
\end{eqnarray}
for $m \geq 1$. So by \eqref{eqn:dertaueiLm},
\begin{eqnarray}
\vl_{1} \left< \tau_{-}^{m}(\vs), \, \ve_{i} \right>
&=&  2(m+3) u_{i}  \left< \tau_{-}^{m}(\vs), \, \ve_{i} \right>
    + \left\{ (m+1)^{2} - \frac{1}{4} \right\}
    \left< \tau_{-}^{m-1}(\vs), \, \ve_{i} \right>   \nonumber \\
&&    + 2(m+1) \sum_{j}  v_{ij} \sqrt{\frac{g_{i}}{g_{j}}}
       \left< \tau_{-}^{m-1}(\vs), \, \ve_{j} \right> \nonumber  \\
&&      + \sum_{j, k}   v_{ij} v_{jk}
       \sqrt{\frac{g_{i}}{g_{k}}}
       \left< \tau_{-}^{m-1}(\vs), \, \ve_{k} \right>
	\label{eqn:L1tauSei}
\end{eqnarray}
for $m \geq 1$.
In particular,
\begin{eqnarray}
\vl_{1} \left< \tau_{-}^{2}(\vs), \, \ve_{i} \right>
&=&  10 u_{i}  \left< \tau_{-}^{2}(\vs), \, \ve_{i} \right>
    + \frac{35}{4}  \sum_{j}  r_{ij} \sqrt{g_{i} \, g_{j}}
    + 6 \sum_{j,k}  v_{ij} r_{jk} \sqrt{g_{i} g_{k}}
\nonumber \\
&&    + \sum_{j, k, l}  v_{ij} v_{jk} r_{kl} \sqrt{g_{i} g_{l}}
    \label{eqn:L1tau2S}
\end{eqnarray}
and
\begin{eqnarray}
\vl_{1} \left( \frac{1} {g_{i}^{2}} \left< \tau_{-}^{3}(\vs), \, \ve_{i} \right>
        \right)
&=&  \frac{63}{4}
    \frac{\left< \tau_{-}^{2}(\vs), \, \ve_{i} \right>}{g_{i}^{2}}
     + 8 \sum_{j}  v_{ij} \frac{1}{\sqrt{g_{i}^{3}g_{j}}}
       \left< \tau_{-}^{2}(\vs), \, \ve_{j} \right>  \nonumber  \\
&&    + \sum_{j, k}   v_{ij} v_{jk}
        \frac{1}{\sqrt{g_{i}^{3}g_{k}}}
       \left< \tau_{-}^{2}(\vs), \, \ve_{k} \right>.
       \label{eqn:L1tau3S}
\end{eqnarray}
The last equation implies that
 $\vl_{1} F_{2}$ does not contain $\left< \tau_{-}^{3}(\vs), \, \ve_{i} \right>$.

Using Lemma~\ref{lem:L1} and equations (\ref{eqn:L1vij}),
(\ref{eqn:L1theta}), (\ref{eqn:L12ndpoleF2}),
(\ref{eqn:L1tau2S}), (\ref{eqn:L1tau3S}), a straightforward computation shows that
\begin{eqnarray}
&& 1152 \, \, \vl_{1} F_{2} = 1152 \, \, \gwiitwo{\vl_{1}}  \nonumber \\
&=& \sum_{i} \frac{\left< \tau_{-}^{2}(\vs), \, \ve_{i} \right>}{g_{i}^{2}}
    \left( 6 + 24 \sum_{j} v_{ij}^{2} \right) \nonumber \\
&& + \sum_{i} \sum_{j \neq i} \theta_{ij} \left\{
    6 \frac{1}{g_{j}}\sqrt{\frac{g_{i}}{g_{j}}}
    + 48 \frac{1}{g_{i}} \sum_{k} v_{ik} v_{jk}
    + 24 \frac{1}{\sqrt{g_{i} g_{j}}} \sum_{k} v_{ik}^{2}
    + 24 \frac{1}{g_{j}}\sqrt{\frac{g_{i}}{g_{j}}} \sum_{k} v_{jk}^{2}
    \right\} \nonumber \\
&& - \sum_{i} 72 r_{ii}^{2} \frac{1}{g_{i}}
   + \sum_{i,j} \left\{ 57 r_{ij}^{2} \frac{1}{g_{i}}
    - 144 r_{ii} r_{ij} v_{ij} \frac{1}{g_{i}} \right\} \nonumber \\
&& + \sum_{i,j,k} \left\{
    \frac{11}{4} r_{ij} r_{ik} \frac{1}{\sqrt{g_{j} g_{k}}}
    + 66 r_{ij} r_{ik} v_{ik} \frac{1}{\sqrt{g_{i} g_{j}}}
    - 36 (r_{ij} v_{ij}) ( r_{ik} v_{ik}) \frac{1}{g_{i}}
    \right. \nonumber \\
&& \hspace{70pt} \left.
    - 288 r_{jk}^{2} v_{ij} v_{ik} \frac{1}{\sqrt{g_{j} g_{k}}}
    + 240 r_{jk}^{2} v_{ij}^{2} \frac{1}{g_{j}}
    \right\}  \nonumber \\
&& + \sum_{i,j,k,l} \left\{
    - 24 r_{jl} r_{kl} v_{ij} v_{ik} \frac{1}{g_{l}}
    + 24 r_{jk} r_{kl} v_{ij}^{2} \frac{1}{\sqrt{g_{j} g_{l}}}
    \right. \nonumber \\
&& \hspace{70pt} \left.
    - 576 r_{jk} (r_{kl} v_{kl}) v_{ij} v_{ik} \frac{1}{g_{k}}
    + 288 r_{jk} (r_{kl} v_{kl}) v_{ij}^{2}
        \frac{1}{\sqrt{g_{j} g_{k}}} \right\} \nonumber \\
&& + \sum_{i,j,k,l,p} \left\{
    - r_{jp} r_{kl} v_{ij} v_{ik} \frac{1}{\sqrt{g_{l} g_{p}}}
    - 24 r_{jp} (r_{kl} v_{kl}) v_{ij} v_{ik}
            \frac{1}{\sqrt{g_{k} g_{p}}}
    \right. \nonumber \\
&& \hspace{70pt} \left.
    - 144 (r_{jp} v_{jp})(r_{kl} v_{kl}) v_{ij} v_{ik}
        \frac{1}{\sqrt{g_{j} g_{k}}} \right\}.
        \label{eqn:L1F2}
\end{eqnarray}

On the other hand, plugging the formulas (\ref{eqn:phi1rot}), (\ref{eqn:phi2iirot}) and
(\ref{eqn:phi2ijrot})
into \eqref{eqn:Virg2L1-3}, then using the fact that
\[ \sum_{i} \sum_{j \neq i} 6 \frac{1}{\sqrt{g_{i} g_{j}}} \theta_{ij}
    = \sum_{i} \sum_{j \neq i} 3 \frac{1}{\sqrt{g_{i} g_{j}}} (\theta_{ij} + \theta_{ji})
    = - \sum_{i} \sum_{j \neq i} 3 \frac{1}{\sqrt{g_{i} g_{j}}} \sum_{k} r_{ik} r_{jk}
 \]
to simplify the resulting expression,
we see that the prediction for the genus-2 $L_{1}$-constraint
is precisely given by \eqref{eqn:L1F2}.
This completes the proof of Theorem~\ref{thm:g2VirSS}.

\appendix
\section{Appendix: Another proof to Theorem \ref{thm:g2VirSS}}

In this appendix, we describe an alternative proof to Theorem \ref{thm:g2VirSS} using
the formula for $\gwiitwo{\vl_{1}}$ given by Theorem~\ref{thm:Lkg2}.
The advantage of this approach is that we do not need to know the precise formula
for expressing $F_{2}$ in terms of rotation coefficients. Therefore
this approach is closer to the treatment of the
genus-1 Virasoro conjecture described in \cite{L4}. Moreover, in this approach,
we are mainly dealing with vector field $T(\bvx)$ which behaves 
much better than $\vl_{1}$. Properties for $T(\bvx)$ derived  here may also
be useful for the study of higher genus Gromov-Witten invariants.

Recall that by Theorem~\ref{thm:Lkg2},
\begin{eqnarray}
 \gwiitwo{\vl_{1}}
     &=&  \frac{1}{2} T(\bvx)
            \left\{ \sum_{i=1}^{N} u_{i}
            B(\ve_{i}, \ve_{i}, \ve_{i}) \right\}
    - \frac{1}{2} \sum_{i=1}^{N} u_{i}^{2}
            B(\ve_{i}, \ve_{i}, \ve_{i})  \nonumber \\
    &&  + A_{1}(\tau_{-}^{2}(\vl_{1}))
        - T(\bvx) \left\{ A_{1} \left( \tau_{-}^{2}(\vl_{0})
            + \frac{3}{2}  \tau_{-}(\vs) \right) \right\}.
    \label{eqn:L1g2}
 \end{eqnarray}
In this section, we will prove that
the right hand side of \eqref{eqn:L1g2} is equal to
the right hand side of \eqref{eqn:L1F2} up to a multiplicative constant 1152.
As observed at the end of last section,
this proves Theorem~\ref{thm:g2VirSS}.

To get rid of $T(\bvx)$ in the
expression of  $\gwiitwo{\vl_{1}}$ in \eqref{eqn:L1g2}, we need
the following properties of this vector field:
\begin{eqnarray}
 T(\bvx) z_{i_{1} \cdots i_{k+2}}
&= &  \sum_{m=2}^{k-1} \,\,\,
       \sum_{1 \leq j_{1} < \cdots < j_{m} \leq k} \,\,\, \sum_{p, q}
    \frac{u_{p}}{g_{q}} \,\,\,
    z_{p i_{j_{1}} \, \cdots \, i_{j_{m}} \, q} \,\,\,
    z_{q \, i_{1} \, \cdots \widehat{i_{j_{1}}} \,
        \cdots \widehat{i_{j_{2}}} \, \cdots \, \cdots\, \cdots
        \widehat{i_{j_{m}}} \, \cdots
        i_{k+2}}  \nonumber \\
&& + \delta_{i_{k+1} i_{k+2}} \sum_{p} u_{p} z_{p \, i_{1} \cdots i_{k+1}}
    - (u_{i_{k+1}} + u_{i_{k+2}}) z_{i_{1} \cdots i_{k+2}}
\end{eqnarray}
and
\begin{eqnarray}
 T(\bvx) \phi_{i_{1} \cdots i_{k}}
&= &  \sum_{m=2}^{k} \,\,\,
 \sum_{1 \leq j_{1} < \cdots < j_{m} \leq k} \,\,\, \sum_{p, q}
    \frac{u_{p}}{g_{q}} \,\,\,
    z_{p i_{j_{1}} \, \cdots \, i_{j_{m}} \, q} \,\,\,
    \phi_{q \, i_{1} \, \cdots \widehat{i_{j_{1}}} \,
        \cdots \widehat{i_{j_{2}}} \, \cdots \, \cdots\, \cdots
        \widehat{i_{j_{m}}} \, \cdots
        i_{k}}  \nonumber \\
&& + \frac{1}{24} \sum_{p, q} \frac{u_{p}}{g_{q}}
     z_{p i_{1} \cdots i_{k} q q}.
\end{eqnarray}
These two equations are obtained from \eqref{eqn:TRRg0}
and (\ref{eqn:TRRg1})
since
\begin{eqnarray*}
 T(\bvx) \gwiig{\vw_{1}, \ldots \vw_{k}}
& =  & \gwiig{T(\bvx) \, \vw_{1} \, \ldots \, \vw_{k}}  \\
&&    + \sum_{j=1}^{k} \gwiig{\vw_{1} \, \ldots \,
(\nabla_{T(\bvx)} \ve_{i})
      \, \ldots \, \vw_{k}},
\end{eqnarray*}
and
 $\nabla_{T(\bvx)} \ve_{i} = - \bvx \qp \ve_{i} = - u_{i} \, \ve_{i}$.
Moreover
\[ T(\bvx) \, g_{i} = - 2 <\bvx, \ve_{i}> = - 2 u_{i} g_{i} \]
for any $i$.

Using the above properties of $T(\bvx)$, we can get rid of this vector field
 in the
right hand side of \eqref{eqn:L1g2}. We separate the contributions
from tensor $A_{1}$ and tensor $B$.
Write
\[ \gwiitwo{\vl_{1}} = L_{A} + L_{B} \]
where
\begin{equation} \label{eqn:LA}
 L_{A}:= A_{1}(\tau_{-}^{2}(\vl_{1}))
        - T(\bvx) \left\{ A_{1} \left( \tau_{-}^{2}(\vl_{0})
                    +\frac{3}{2} \tau_{-}(\vs)  \right) \right\}
\end{equation}
 is the contribution from the tensor $A_{1}$, and
\begin{equation} \label{eqn:LB}
 L_{B} := \frac{1}{2} \sum_{i}
               \left(u_{i} T(\bvx) B(\ve_{i}, \ve_{i}, \ve_{i})
                - u_{i}^{2} B(\ve_{i}, \ve_{i}, \ve_{i})  \right)
\end{equation}
is the contribution from the tensor $B$.
We have
\begin{eqnarray*}
2 L_{B} &=&
 \sum_{i,j,k} \frac{u_{i} u_{j}}{g_{j} g_{k}} \, \, \frac{1}{240} z_{ijjjjkk}
    +  \sum_{i,j,k, p} \frac{u_{i} u_{p}}{g_{j} g_{k}}
    \left( - \frac{1}{480} z_{iijjkkp}
    - \frac{1}{2880} z_{iiijkkp} \right) \\
&& + \sum_{i,j,k, p, q} \frac{u_{i} u_{p}}{g_{j} g_{k} g_{q}}
    \left( \frac{1}{2880} z_{iiijjk} z_{kpqq}
    - \frac{1}{480} z_{iiij} z_{jkkpqq}
    + \frac{1}{320} z_{ijjk} z_{iikpqq}  \right. \\
&& \left. \hspace{80pt}    - \frac{1}{80} z_{iijk} z_{ijkpqq}
    + \frac{1}{240} z_{iiijk} z_{jkpqq}
    - \frac{1}{320} z_{iijjk} z_{ikpqq}
    \right) \\
&& + \sum_{i} \phi_{i} \left\{
     \sum_{j,k} \left(
         \frac{1}{10} \frac{u_{i} u_{j}}{g_{i} g_{k}} z_{iiijkk}
     + \frac{1}{10} \frac{u_{j} u_{k}}{g_{i} g_{k}} z_{ijkkkk}
     + \frac{1}{120} \frac{u_{k}(2u_{i} + 2 u_{j} - 3 u_{k})}{g_{i} g_{j}}
                  z_{ijjkkk} \right) \right. \\
&&  \hspace{35pt}   - \sum_{j,k,p} \frac{u_{k} u_{p}}{g_{i} g_{j}} \left(
               \frac{11}{120} z_{ijjkkp} + \frac{1}{120} z_{ijkkkp}
           \right) \\
&&  \hspace{35pt}
           + \sum_{j,k,p,q} \frac{u_{p} u_{q}}{g_{i} g_{j} g_{k}} \left(
               \frac{1}{40} z_{ijqqq} z_{jkkp}
           - \frac{11}{120} z_{ijkkp} z_{jqqq}
           -  \frac{1}{24} z_{ikqq} z_{jjkpq} \right. \\
&& \hspace{120pt}
           + \frac{7}{60} z_{ijkp} z_{jkqqq}
           + \frac{1}{60} z_{ijkqq} z_{jkpq}
           - \frac{1}{15} z_{ikpq} z_{jjkqq}  \\
&& \hspace{120pt} \left. \left.
           -  \frac{17}{60} z_{ijkpq} z_{jkqq}
           + \frac{3}{40} z_{ikpqq} z_{jjkq}
           \right)
      \right\} \\
&& - \sum_{i} \phi_{i i} \sum_{j, k}
       \frac{3}{10}  \frac{u_{i} u_{j}}{g_{i} g_{k}} z_{iijkk} \\
&&      + \sum_{i, j} \phi_{i j} \left\{
       \sum_{k} \left(
           \frac{2}{5} \frac{u_{i} u_{k}}{g_{i} g_{j}} z_{iiijk}
           + \frac{1}{10}
           \frac{u_{k}(2u_{i} + 2 u_{j} - 3 u_{k})}{g_{i} g_{j}}
                  z_{ijkkk} \right. \right. \\
&& \hspace{100pt} \left.   - \frac{3}{40}
           \frac{u_{i}(2u_{j} + 2 u_{k} - 3 u_{i})}{g_{j} g_{k}}
                  z_{iijkk} \right) \\
&&  \hspace{30pt}
           - \sum_{k, p} \left(
           \frac{1}{40} \frac{u_{i} u_{p}}{g_{j} g_{k}}
                (z_{ijkkp} + z_{iijkp})
           + \frac{1}{120}
           \frac{u_{k}  u_{p}}{g_{i} g_{j}}
                  z_{ikkkp} \right)  \\
&& \hspace{30pt}
       + \sum_{k, p, q} \left(
            \frac{1}{40} \frac{u_{p} u_{q}}{g_{i} g_{j} g_{k}}
                \left( - 4 z_{ikkp} z_{jqqq} - 16 z_{ikpq} z_{jkqq}
            + 3 z_{ipqq} z_{jkkq}
            \right)    \right. \\
&& \hspace{80pt} \left. \left.
           + \frac{1}{40}
           \frac{u_{i}  u_{p}}{g_{j} g_{k} g_{q}}
                  \left(- 5 z_{iijk} z_{kpqq}
               - 18 z_{iikq} z_{jkpq} + 6 z_{ikqq} z_{ijkp}
               \right)
              \right)
      \right\} \\
&& + \sum_{i} \phi_{i i i} \sum_{j, k}
       \frac{1}{10}  \frac{u_{i} u_{j}}{g_{i} g_{k}} z_{ijkk} \\
&&  + \sum_{i, j} \phi_{i i j} \left\{
       \sum_{k} \left(
           \frac{3}{5} \frac{u_{i} u_{k}}{g_{i} g_{j}} z_{iijk}
           - \frac{1}{20}
           \frac{u_{k}(2u_{i} + 2 u_{j} - 3 u_{k})}{g_{i} g_{j}}
                  z_{jkkk} \right. \right. \\
&& \hspace{100pt} \left.   + \frac{3}{40}
           \frac{u_{i}(2u_{j} - u_{i})}{g_{j} g_{k}}
                  z_{ijkk} \right) \\
&&  \hspace{30pt} \left.
           + \sum_{k, p} \left(
           - \frac{1}{10} \left( \frac{u_{i} u_{p}}{g_{j} g_{k}}
                  + \frac{u_{k} u_{p}}{g_{i} g_{j}} \right) z_{jkkp}
           + \frac{1}{20}
           \frac{u_{i}  u_{p}}{g_{j} g_{k}}
                  z_{ijkp} \right) \right\} \\
&& + \sum_{i, j, k} \phi_{i j k} \left\{
            -  \frac{3}{10}
           \frac{u_{i}(2u_{j} + 2 u_{k} - 3 u_{i})}{g_{j} g_{k}}
                  z_{iijk}
            - \sum_{p} \frac{u_{i} u_{p}}{g_{j} g_{k}}
                \left( \frac{1}{40} z_{iikp}
                     + \frac{1}{2} z_{ijkp} \right) \right\} \\
&& + \sum_{i} \phi_{iiii} \, \cdot \,
     \frac{1}{10} \frac{u_{i}^{2}}{g_{i}}  -
     \sum_{i, j} \left(\frac{1}{120}  \, \phi_{iiij}
             +  \frac{1}{20} \, \phi_{iijj} \right) \,
         \frac{u_{i} (2 u_{j} - u_{i})}{g_{j}} \\
&& + \sum_{i, j} \phi_{i} \phi_{j} \left\{
       \sum_{k} \left(
           \frac{12}{5} \frac{u_{i} u_{k}}{g_{i} g_{j}} z_{iiijk}
           + \frac{1}{5}
           \frac{u_{k}(2u_{i} + 2 u_{j} - 3 u_{k})}{g_{i} g_{j}}
                  z_{ijkkk} \right) \right. \\
&&  \hspace{50pt} \left.
           - \sum_{k, p}  \frac{u_{k}  u_{p}}{g_{i} g_{j}}
                  z_{ijkkp}
       - \sum_{k, p, q} \frac{u_{p} u_{q}}{g_{i} g_{j} g_{k}}  \left(
            \frac{4}{5}
                 z_{ikpq} z_{jkqq} +  z_{ijkp} z_{kqqq}
              \right)
      \right\} \\
&& - \sum_{i, j} \phi_{i} \phi_{jj} \sum_{k}
       \frac{36}{5}  \frac{u_{j} u_{k}}{g_{i} g_{j}} z_{ijjk} \\
&&  + \sum_{i, j} \phi_{i} \phi_{i j} \left\{
       \sum_{k} \left(
           \frac{36}{5} \frac{u_{i} u_{k}}{g_{i} g_{j}} z_{iijk}
           - \frac{6}{5}
           \frac{u_{k}(2u_{i} + 2 u_{j} - 3 u_{k})}{g_{i} g_{j}}
                  z_{jkkk}  \right) \right. \\
&&  \hspace{65pt} \left.
           - \sum_{k, p} \frac{12}{5} \frac{u_{k} u_{p}}{g_{i} g_{j}}
                   z_{jkkp}  \right\} \\
&& - \sum_{i, j, k} \phi_{i} \phi_{j k} \,\, \cdot \,\,
             \frac{6}{5}
           \frac{u_{j}(2u_{i} + 2 u_{k} - 3 u_{j})}{g_{i} g_{k}}
                  z_{ijjk} \\
&& + \sum_{i} \phi_{i} \phi_{iii} \, \cdot \,
         \frac{12}{5} \,  \frac{u_{i}^{2}}{g_{i}}
    - \sum_{i, j} \phi_{i} \phi_{ijj}  \, \cdot \,
    \frac{6}{5} \, \frac{u_{j}(2u_{i}-u_{j})}{g_{i}}  \\
&&
    - \sum_{i} \phi_{ii}^{2} \, \cdot \,
         \frac{18}{5} \,  \frac{u_{i}^{2}}{g_{i}}
    + \sum_{i, j} \phi_{ij}^{2}  \, \cdot \,
    \frac{9}{5} \, \frac{u_{i}(2u_{j}-u_{i})}{g_{j}}.
\end{eqnarray*}

To compute the contribution from tensor $A_{1}$,
we also need to understand the action of $T(\bvx)$ on coefficients
of $\tau_{-}^{k}(\vl_{0})$ and $\tau_{-}^{k}(\vs)$ along idempotents.
Since
\[ \nabla_{T(\vw)} \tau_{-}^{k}(\vs)
= \tau_{-}^{k}(\nabla_{T(\vw)} \vs) = - \tau_{-}^{k+1}(T(\vw))
= - \tau_{-}^{k}(\vw),  \]
by \eqref{eqn:dertauei},
\[ T(\vw) <\tau_{-}^{k}(\vs), \, \ve_{i}>
    = - <\tau_{-}^{k}(\vw), \, \ve_{i}>
        - <\tau_{-}^{k}(\vs), \, \vw \qp \ve_{i}>  \]
for any vector field $\vw$ and $k \geq 0$. In particular, for
$k \geq 1$,
\begin{equation} \label{eqn:TXtauS}
 T(\bvx) <\tau_{-}^{k}(\vs), \, \ve_{i}> = - u_{i}
    <\tau_{-}^{k}(\vs), \, \ve_{i}>.
\end{equation}
Comparing to \eqref{eqn:L1tauSei}, we see that $T(\bvx)$ 
behaves much better than $\vl_{1}$.
Similarly, since $\nabla_{T(\vw)} \vl_{0} = R(\vw)$
(cf. \cite[Equation (43)]{L2}),
\[ T(\vw) <\tau_{-}^{k}(\vl_{0}), \, \ve_{i}>
    =  <\tau_{-}^{k}R(\vw), \, \ve_{i}>
        - <\tau_{-}^{k}(\vl_{0}), \, \vw \qp \ve_{i}>  \]
for any vector field $\vw$.
By \cite[Lemma 3.11]{L2},
for primary vector field $\bvw$ and $k \geq 2$,
\[ \tau_{-}^{k}R(\bvw)
    = \tau_{-}^{k-1} \left(R \tau_{-}(\bvw) + \vg * \bvw + \bvw \right)
    = 0. \]
So for $k \geq 2$,
\begin{equation} \label{eqn:TXtauL0}
 T(\bvx) <\tau_{-}^{k}(\vl_{0}), \, \ve_{i}> = - u_{i}
    <\tau_{-}^{k}(\vl_{0}), \, \ve_{i}>.
\end{equation}
Therefore for
\[ \vw =  \tau_{-}^{k}(\vs) \,\,\, {\rm or} \,\,\,
             \tau_{-}^{k+1}(\vl_{0}) \,\,\, {\rm with} \,\,\,
         k \geq 1,\]
we can use \eqref{eqn:TXtauS} and \eqref{eqn:TXtauL0}
 to compute $T(\bvx) A_{1}(\vw)$
and obtain
\begin{eqnarray*}
T(\bvx) A_{1}(\vw)
&=& \sum_{i}  \frac{u_{i} \left< \vw, \, \ve_{i} \right>}{g_{i}^{2}}
            \left( \frac{21}{10} \phi_{i}^{2}
                + \frac{3}{10} \phi_{ii} \right)
    - \sum_{i, j}  \frac{1}{240}
                   \frac{\left< \vw, \, \ve_{i} \right>}{g_{i}}
           \frac{u_{i}+2 u_{j}}{g_{j}} \phi_{ij}
                 \\
&&  + \sum_{i}  \phi_{i} \left\{
        \frac{3}{20}
           \frac{u_{i} \left< \tau_{-}(\vw), \, \ve_{i} \right>}{g_{i}^{2}}
        + \sum_{j, k, p} \frac{1}{20}
                  \frac{\left< \vw, \, \ve_{k} \right>}{g_{k}}
               \frac{u_{p}}{g_{i} g_{j}} z_{ijkp}
            \right.  \\
&& \hspace{50pt}
            + \sum_{j, k} \left( \left(\frac{7}{120}
                  \frac{\left< \vw, \, \ve_{i} \right>}{g_{i}}
              + \frac{1}{10}
                  \frac{\left< \vw, \, \ve_{k} \right>}{g_{k}}
              \right)
            \frac{u_{j}}{g_{i} g_{k}} z_{ijkk}    \right.  \\
&& \hspace{90pt} \left. \left.
        + \frac{13}{240}  \frac{\left< \vw, \, \ve_{k} \right>}{g_{k}}
               \frac{2u_{i} + u_{k}}{g_{i} g_{j}} z_{ijjk}
               \right)
             \right\} \\
&& + \sum_{i} \frac{\left< \vw, \, \ve_{i} \right>}{g_{i}}
         \left\{ \sum_{j,k} \left(
                  \frac{1}{240} \frac{u_{j}}{g_{i} g_{k}} z_{iijkk}
                + \frac{1}{960} \frac{u_{i} + 2u_{j}}{g_{j} g_{k}}
                     z_{ijjkk} \right)  \right.  \\
&& \hspace{80pt} \left.
        + \sum_{j,k,p}  \frac{1}{1152}
                 \frac{u_{p}}{g_{j} g_{k}} z_{ijkkp}
                \right. \\
&& \hspace{80pt} \left.
    + \sum_{j, k, p, q}  \frac{1}{480} \frac{u_{p}}{g_{j} g_{k} g_{q}} \left(
            \frac{19}{12}   z_{ijjk} z_{kpqq}
               +  z_{ijpq} z_{jkkq} \right) \right\} \\
&& + \sum_{i}  \frac{\left< \tau_{-}(\vw), \, \ve_{i} \right>}{g_{i}}
           \left\{ \sum_{j}  \frac{1}{160} \frac{u_{i}}{g_{i} g_{j}} z_{iijj}
        + \sum_{j,k,p}  \frac{1}{1152}
                 \frac{u_{p}}{g_{j} g_{k}} z_{ijkp}
                \right. \\
&& \hspace{65pt} \left.
    + \sum_{j,k} \left(
                  \frac{1}{480} \frac{u_{j}}{g_{i} g_{k}} z_{ijkk}
                + \frac{1}{1152} \frac{u_{i} + 2u_{j}}{g_{j} g_{k}}
                     z_{ijkk}
            + \frac{1}{480} \frac{u_{k}}{g_{i} g_{j}} z_{iijk}
             \right)
             \right\} \\
&&    + \sum_{i}  \frac{\left< \tau_{-}^{2}(\vw), \, \ve_{i} \right>}{g_{i}}
             \,\, \cdot \,\,  \frac{1}{384}
               \frac{u_{i}}{g_{i}}.
\end{eqnarray*}
Define
\begin{eqnarray}
 c_{i j; k} &:= &  \frac{1}{g_{i}} \left\{
               \left< \tau_{-}^{k}(\vl_{1}), \,  \ve_{i} \right>
           - (u_{i} + 2 u_{j})
           \left< \tau_{-}^{k}(\vl_{0}) + \frac{3}{2} \tau_{-}^{k-1}(\vs),
               \,  \ve_{i} \right> \right\} \label{eqn:cijk}
\end{eqnarray}
and
\begin{eqnarray}
 d_{i; k} &:= &  \frac{1}{g_{i}}
           \left< \tau_{-}^{k}(\vl_{0}) + \frac{3}{2} \tau_{-}^{k-1}(\vs),
               \,  \ve_{i} \right>. \label{eqn:dij}
\end{eqnarray}
By equations (\ref{eqn:L0idem}) and (\ref{eqn:taukL1ei}),
\begin{eqnarray*}
 c_{i j; k}
&=& 2 u_{i} u_{j} \frac{ \left< \tau_{-}^{k}(\vs), \,  \ve_{i} \right>}{g_{i}}
    - \left\{(k+2) u_{i} + 2 (1-k) u_{j} \right\}
    \frac{ \left< \tau_{-}^{k-1}(\vs), \,  \ve_{i} \right>}{g_{i}} \\
&&    + \sum_{p} (u_{p} - 2 u_{j}) (u_{i} - u_{p}) r_{ip}
              \frac{1}{\sqrt{g_{i} g_{p}}}
          \left< \tau_{-}^{k-1}(\vs), \,  \ve_{p} \right> \\
&&    + (\frac{1}{4} - k^{2})
      \frac{ \left< \tau_{-}^{k-2}(\vs), \,  \ve_{i} \right>}{g_{i}}
    + 2k \sum_{p}  (u_{i} - u_{p}) r_{ip}
              \frac{1}{\sqrt{g_{i} g_{p}}}
          \left< \tau_{-}^{k-2}(\vs), \,  \ve_{p} \right> \\
&&    - \sum_{p, q} (u_{i} -  u_{p}) (u_{p} - u_{q}) r_{ip} r_{pq}
              \frac{1}{\sqrt{g_{i} g_{q}}}
          \left< \tau_{-}^{k-2}(\vs), \,  \ve_{q} \right>
\end{eqnarray*}
and
\begin{eqnarray*}
 d_{i; k}
&=& - u_{i}  \frac{ \left< \tau_{-}^{k}(\vs), \,  \ve_{i} \right>}{g_{i}}
    + (1-k)  \frac{ \left< \tau_{-}^{k-1}(\vs), \,  \ve_{i} \right>}{g_{i}}
    -   \sum_{p} (u_{p} - u_{i})  r_{ip}
              \frac{1}{\sqrt{g_{i} g_{p}}}
          \left< \tau_{-}^{k-1}(\vs), \,  \ve_{p} \right>.
\end{eqnarray*}

The contribution to $\gwiitwo{\vl_{1}}$ from tensor $A_{1}$ is
\begin{eqnarray*}
L_{A}
&=& \sum_{i} \frac{7}{10} \frac{1}{g_{i}} c_{ii;2} \phi_{i}^{2}
    +  \sum_{i} \frac{1}{10} \frac{1}{g_{i}} c_{ii;2} \phi_{ii}
    - \sum_{i, j}  \frac{1}{240} \frac{1}{g_{j}} c_{ij;2} \phi_{ij} \\
&&    + \sum_{i} \phi_{i} \left\{
    \frac{1}{20} \frac{1}{g_{i}} c_{ii;3}
    + \sum_{j,k} \left(
                 \frac{13}{240} \frac{1}{g_{i} g_{j}} c_{ki;2} z_{ijjk}
                -  \frac{7}{120} \frac{u_{j}}{g_{i} g_{k}} d_{i;2} z_{ijkk}
        -  \frac{1}{10} \frac{u_{j}}{g_{i} g_{k}} d_{k;2} z_{ijkk}
        \right) \right. \\
&& \hspace{60pt} \left.
   - \sum_{j,k,p} \frac{1}{20} \frac{u_{p}}{g_{i} g_{j}} d_{k;2} z_{ijkp}
    \right\} \\
&& + \sum_{i}  \frac{1}{1152} \, \frac{1}{g_{i}} c_{ii;4}
    + \sum_{i, j} \frac{1}{480} \frac{1}{g_{i} g_{j}} c_{ii;3} z_{iijj} \\
&&     +  \sum_{i, j, k} \frac{1}{1152} \,
              \frac{1}{g_{j} g_{k}} c_{ij;3} z_{ijkk}
    +  \sum_{i, j, k} \frac{1}{960} \frac{1}{g_{j} g_{k}} c_{ij;2} z_{ijjkk} \\
&&  - \sum_{i} d_{i;3} \left\{ \sum_{j, k} \frac{1}{480}
              \left( \frac{u_{j}}{g_{i} g_{k}} z_{ijkk} +
                 \frac{u_{k}}{g_{i} g_{j}} z_{iijk} \right)
         + \sum_{j, k, p} \frac{1}{1152} \,
                   \frac{u_{p}}{g_{j} g_{k}} z_{ijkp} \right\} \\
&&  - \sum_{i} d_{i;2} \left\{
           \sum_{j, k} \frac{1}{240} \frac{u_{j}}{g_{i} g_{k}}
                    z_{iijkk}
            + \sum_{j, k, p} \frac{1}{1152} \,
                   \frac{u_{p}}{g_{j} g_{k}} z_{ijkkp}
           \right. \\
&& \hspace{80pt} \left.
           + \sum_{j, k, p, q} \frac{1}{480}
                      \frac{u_{p}}{g_{j} g_{k} g_{q}}
           \left(\frac{19}{12}  z_{ijjk} z_{kpqq}
           + z_{ijpq} z_{jkkq} \right)
\right\}.
\end{eqnarray*}

In section \ref{sec:RotCoeff}, we have described how to represent
genus-0 and genus-1 functions $z_{i_{1} \cdots i_{k}}$
and $\phi_{i_{1} \cdots i_{k}}$ in terms of rotation coefficients.
For our purpose, we only need $z_{i_{1} \cdots i_{k}}$ for
$4 \leq k \leq 7$ and $\phi_{i_{1} \cdots i_{k}}$ for $1 \leq k
\leq 4$. Using these formulas, we can express $L_{A} + L_{B}$ in terms of
functions $u_{i}$, $g_{i}$, $r_{ij}$, and $ <\tau_{-}^{k}(\vs),
\, \ve_{i}>$ with $k \geq 2$. After lengthy but
straightforward computations, we can check that $L_{A} + L_{B}$
is equal to the right hand side of \eqref{eqn:L1F2} up to a multiplicative
constant 1152, and
therefore  Theorem~\ref{thm:g2VirSS} is proved. It is interesting
to observe what happens to terms $ <\tau_{-}^{k}(\vs), \,
\ve_{i}>$ in the computation. Both $L_{A}$ and $L_{B}$ contains
terms $ <\tau_{-}^{k}(\vs), \, \ve_{i}>$ with $2 \leq k \leq 4$.
But coefficients of $ <\tau_{-}^{3}(\vs), \, \ve_{i}>$ (as well as
$ <\tau_{-}^{4}(\vs), \, \ve_{i}>$)  from $L_{A}$ and $L_{B}$ are
opposite to each other, therefore are cancelled in $L_{A} +
L_{B}$. However $L_{A} + L_{B}$ contains $ <\tau_{-}^{2}(\vs), \,
\ve_{i}>$, which exactly match the corresponding terms in
 the right hand side of \eqref{eqn:L1F2}. In a contrast, the formula for $F_{2}$
in Theorem~\ref{thm:F2rot} contains
 $ <\tau_{-}^{3}(\vs), \, \ve_{i}>$. The action of $\vl_{1}$ transforms this term to
expressions only involve $ <\tau_{-}^{2}(\vs), \, \ve_{i}>$ as indicated in
\eqref{eqn:L1tau3S}.


\vspace{30pt}
\noindent
Department of Mathematics  \\
University of Notre Dame \\
Notre Dame,  IN  46556, USA \\

\vspace{10pt} \noindent E-mail address: {\it xliu3@nd.edu}


\begin{thebibliography}{399}

\bibitem[BP]{BP} Belorousski, P. and Pandharipande, R.,
        {\it A descendent relation in genus 2},
    Ann. Scuola Norm. Sup. Pisa Cl. Sci. (4) 29 (2000) 171-191.
\bibitem[CK]{CK} Cox, D. and Katz, S.,
        {\it Mirror symmetry and algebraic geometry},
        Providence, R.I. AMS, 1999.
\bibitem[D]{D} B. Dubrovin,
    {\it Geometry of 2D topological field theories},
    Integrable Systems and Quantum Groups,
    Springer Lectures Notes in Math. 1620 (1996), 120-348.
\bibitem[DZ1]{DZ1} Dubrovin, B., Zhang, Y.,
           {\it Bihamiltonian hierarchies in 2D topological field
                theory at one-loop approximation},
            Comm. Math. Phys. 198 (1998), no.2, 311-361.
\bibitem[DZ2]{DZ2} Dubrovin, B., Zhang, Y.,
           {\it Frobenius manifolds and Virasoro constraints}, \\
           Selecta Math. (N.S.) 5 (1999) 423-466.
\bibitem[DZ3]{DZ3} Dubrovin, B., Zhang, Y.,
           {\it Normal forms of hierarchies of integrable PDEs,
           Frobenius manifolds and Gromov-Witten invariants}, \\
           math.DG/0108160.
\bibitem[EHX]{EHX} Eguchi, T., Hori, K., and Xiong, C.,
        {\it Quantum Cohomology and Virasoro Algebra},
        Phys. Lett. B402 (1997) 71-80.
\bibitem[Ge1]{Ge1} Getzler, E.,
        {\it Topological recursion relations in genus 2},
       Integrable systems and algebraic geometry (Kobe/kyoto, 1997)
    73-106.
\bibitem[Ge2]{Ge2} Getzler, E.,
          {\it The Virasoro conjecture for Gromov-Witten invariants}, \\
          ({\textsf math.AG/9812026})
\bibitem[Gi1]{Gi1} Givental, A.,
        {\it Semisimple Frobenius structures at higher genus},
        Intern. Math. Res. Notices, 23 (2001), 1265-1286.
\bibitem[Gi2]{Gi2} Givental, A.,
        {\it Gromov-Witten invariants and quantization of quadratic
        hamiltonians},
        Moscow Mathematical Journal, v.1, no. 4 (2001), 551-568.
\bibitem[Gi3]{Gi3} Givental, A.,
        {\it Symplectic geometry of Frobenius structures},
        math.AG/0305409.
\bibitem[Ko]{Kon} Kontsevich, M.,
        {\it Intersection theory on the moduli space of curves and
                the matrix airy function}, Comm. Math. Phys., 147, 1-23 (1992).
\bibitem[LiT]{LiT} Li, J. and Tian, G.,
        {\it Virtual moduli cycles and Gromov-Witten invariants of
                general symplectic manifolds},
        Topics in symplectic 4-manifolds (Irvine, CA, 1996), 47-83.
\bibitem[L1]{L1} X. Liu,
    {\it Elliptic Gromov-Witten invariants and Virasoro conjecture},
    Comm. Math. Phys. 216 (2001), 705-728.
\bibitem[L2]{L2} X. Liu,
    {\it Quantum product on the big phase space and Virasoro conjecture},
    Advances in Mathematics 169 (2002), 313-375.
\bibitem[L3]{L3} Liu, X.,
    {\it Quantum product, topological recursion relations, and the
    Virasoro conjecture}, to appear in Preceedings of Mathematical Society
    of Japan - 9th International Research Institute on ``Integrable Systems
    in Differential Geometry'' in 2000, Tokyo, Japan.
\bibitem[L4]{L4} X. Liu,
    {\it Idempotents on the big phase space},
    math.DG/0310409.
\bibitem[LT]{LT} X. Liu and G. Tian,
        {\it Virasoro constraints for quantum cohomology},  \\
        J. Diff. Geom. 50 (1998), 537 - 591.
\bibitem[RT]{RT} Ruan, Y. and  Tian, G.,
        {\it Higher genus symplectic invariants and sigma models coupled
         with gravity}, Invent. Math. 130 (1997), 455-516.
\bibitem[W1]{W1} Witten, E.,
        {\it Two dimensional gravity and intersection theory on
                Moduli space},
        Surveys in Diff. Geom., 1 (1991), 243-310.
\bibitem[W2]{W2} Witten, E.,
          {\it On the Kontsevich model and other models of two dimensional
           gravity},
          in "Proceedings of the XXth international conference on differential
          geometric methods in theoretical physics (New York, 1991)",
          World Sci. Publishing, River Edge, NJ, 1992, pp. 176-216.
\end{thebibliography}
\end{document}